\documentclass{ip-journal}

\newtheorem{theorem}{Theorem}[section]
\newtheorem{corollary}[theorem]{Corollary}
\newtheorem{definition}[theorem]{Definition}
\newtheorem{lemma}[theorem]{Lemma}
\renewcommand{\qed}{\mbox{}\hfill $\square$}

\newcommand{\T}{\top}
\def\qn#1{\big|\!\big|\!\big|\,#1\,\big|\!\big|\!\big|}
\def\mc{\mathcal}
\def\com#1{\quad\text{#1}\quad}
\def\wint{\star\!\!\int}
\newcommand{\wt}{\widetilde}
\def\del{\partial}
\newcommand{\DD}{\mathbb{D}}
\newcommand{\D}{\mbox{D}}
\newcommand{\BB}{\mathfrak{B}}
\newcommand{\RR}{\mathbb{R}}
\newcommand{\eps}{\varepsilon}
\newcommand{\Var}{\mathbb{V}}
\newcommand{\<}{\big\langle}
\renewcommand{\>}{\big\rangle}
\DeclareMathOperator*{\esssup}{ess\,sup}

\title[Weak* solutions I]{Weak* solutions I: A new perspective on solutions to
 systems of conservation laws}

\author{Alexey Miroshnikov \and   Robin Young }

\address{
Department of Biostatistics and Epidemiology, Univ of
Massachusetts Amherst \\ current address: Department of Mathematics,
UCLA \\  Department of Mathematics and Statistics,
Univ of Massachusetts Amherst}

\email{amiroshn@gmail.com,  young@math.umass.edu}



\begin{document}

\begin{abstract}
  We introduce a new notion of solution, which we call weak*
  solutions, for systems of conservation laws.  These solutions can be
  used to handle singular situations that standard weak solutions
  cannot, such as vacuums in Lagrangian gas dynamics or cavities in
  elasticity.  Our framework allows us to treat the systems as ODEs in
  Banach space.  Starting with the observation that solutions act
  linearly on test functions $\alpha\in X$, we require solutions to
  take values in the dual space $X^*$ of $X$.  Moreover, we weaken the
  usual requirement of measurability of solutions.  In order to do
  this, we develop the calculus of the Gelfand integral, which is
  appropriate for weak* measurable functions.  We then use the Gelfand
  calculus to define weak* solutions, and show that they are stronger
  than the usual notion of weak solution, although for $BV$ solutions
  the notions are equivalent.  It is expected that these solutions
  will also shed light on vexing issues of ill-posedness for
  multi-dimensional systems.
\end{abstract}

\maketitle

\section{Introduction}

Hyperbolic systems of conservation laws are fundamental in the study
of fluids and continuum dynamics, expressing basic physical properties
such as conservation of mass and momentum.  In one dimension, these
take the form
\begin{equation}
  u_t + F(u)_x = 0, \qquad u(x,0) = u^0(x),
\label{cl}
\end{equation}
where $u\in \mc U\subset\RR^n$, $F:\mc U\to\RR^n$, being derived from
the integral form
\[
  \frac{d}{dt}\int_{(a,b)} u\; dx = - \int_{\partial(a,b)}F\cdot n\;ds
  = - F(u)\Big|_a^b.
\]
Here $u$ is the vector of conserved quantities, and $F$ is the
corresponding flux.  Mathematically, these equations present
particular challenges because of the ubiquity of shock waves, at which
continuum fields such as pressure are discontinuous, so that global
classical solutions generally do not exist.

The study of weak solutions for general systems was initiated by Peter
Lax in the 1950's~\cite{Lax}, and followed in the 1960's by
Glimm's celebrated proof of the global existence of weak solutions,
provided the initial data has sufficiently small total
variation~\cite{G}.  A weak solution is a locally bounded, locally
integrable function $u$ satisfying
\begin{equation}
  \label{weak}
\begin{aligned}
  \int_0^T\int_{\RR} \Big( u(x,t)^{\T} \del_t\varphi(x,t)&
  + F(u(x,t))^{\T} \del_x\varphi(x,t) \Big)\; dx\;dt\\
   &{} + \int_{\RR}  u^0(x)^{\T} \varphi(0,x) \; dx = 0,
\end{aligned}
\end{equation}
for all smooth test functions $\varphi$~\cite{Lax,Dafermos}.  This
must be augmented by an \emph{entropy condition}, which serves to
restrict weak solutions to those which contain only admissible shocks.
There is now a mature and complete theory of entropy weak solutions in
one space dimension, culminating in Bressan and collaborators' proofs
of uniqueness and continuous dependence of entropy weak solutions,
both via approximation schemes and via vanishing
viscosity~\cite{Bressan, BB}.  All of these results depend on Glimm's
global $BV$ estimates, which in turn require the initial data to have
small amplitude and bounded variation~\cite{Yth}.

If the data is large, even for physical systems, the class of weak
solutions is not always sufficient to cope with the variety of
solutions.  In particular, the use of a Lagrangian frame, which is
often computationally convenient, presents problems near vacuums,
cavities and fractures in the medium, where they are represented by
Dirac masses~\cite{YpRP2, Ycoll}.  There have been several varied
attempts to extend the class of weak solutions in these and similar
contexts~\cite{KK2, DanShel, Sever, GT}.  This problem is worse in
higher dimensions, for which there are several instances of solutions
exhibiting instability and nonuniqueness, independent of entropy
considerations~\cite{Rauch, Scheffer, Shn, DeLellis}.

Our initial motivation was to allow Dirac masses, representing
vacuums, which are discontinuities in the medium, to appear in
solutions of the Euler equations of gas dynamics, in a Lagrangian
frame.  In particular, once the use of Dirac masses is rigorously
justified, their use simplifies the manipulation of solutions in
calculations.  In addition, we would like to be able to rigorously
describe solutions as the solutions of an ODE in a Banach space, thus
correctly regarding the conservation law as an evolution equation.
We would also prefer a definition that more closely resembles the
integral form of the balance law than does \eqref{weak}.

The pairing \eqref{weak} of the weak solution $u(x,t)$ with arbitrary
smooth test functions $\varphi \in C_c^{\infty}(\RR \times (0,T))$,
together with the idea of studying the evolution
$t \to u(t)=u(\cdot, t)$, motivate the following observations.
\begin{itemize}
\item $u(t)$, the value of a solution $u$, should be a spatial
  function which acts linearly on spatial test functions
  $x\to\alpha(x)$, that live in a certain Banach space $X$. Thus,
  $u(t)$ should take values in the dual space, $u(t) \in X^*$.
\item Rather than requiring $t \to u(t)\in X^*$ to be (strongly)
  Bochner integrable, it is sufficient to require the weaker condition
  that the scalar function $t \to \<u(t),\alpha\>$ be integrable for
  each test function $\alpha\in X$.
\item To make sense of the nonlinear flux $F(u)$ and its derivative in
  the space $X^*$, one needs to properly define the flux as a mapping
  $F:X^* \to X^*$, so $F$ must be extended to a map of Banach spaces.
\item To obtain consistency and be able to use the calculus of
  distributions, we also should require that a union of sets
  $C_c^{\infty}(\Omega)$ are dense in $X$ for appropriate subsets
  $\Omega \subseteq \RR$ .
\end{itemize} 

In this paper, we develop the point of view that $u=u(t)\in X^*$ is a
map on $[0,T]$ with values in a Banach space $X^*$, which solves the
evolution equation
\begin{equation}\label{BODE}
  u'(t) + \DD_x F(u(t)) = 0 \com{in} X^*,
\end{equation}
in a natural sense, for appropriately interpreted time derivative $u'$
and flux $F(u)$.  That is, we view the problem as a time evolution,
treating space using the relation between the test function space
$X$ and its dual $X^*$, while for time we use the usual test functions
$C_c^\infty(0,T)$.  We note that weak solutions, defined by
\eqref{weak}, make no distinction between time and space.

In light of our second observation above, we regard the assumption of
strong measurability, in which $u:[0,T]\to X^*$ is measurable, as
being too strong for our purposes.  Instead, we prefer to use the
concept of weak* measurability, so that the scalar function
$\<u,\alpha\>:[0,t]\to \RR$ is measurable for each $\alpha\in X$.  The
appropriate notion of integration in this case is the Gelfand
integral, denoted $\wint_{E} u(t) \, dt$, with $E \in \BB([0,T])$
measurable, which is essentially defined by the condition
\[
  \wint_{E} u(t)\;dt \in X^*,   \com{such that} 
  \Big\langle \wint_{E} u(t)\;dt,\alpha\Big\rangle
  = \int_{E} \<u(t),\alpha\>\;dt,
\]
for all spatial test functions $\alpha\in X$.  In fact, in order to
make sense of the solution $u(t)$ as a measure, we are required to use
the Gelfand integral rather than the Bochner integral.

Although the Gelfand integral goes back to the 1930's, to the authors'
knowledge it has not been extensively studied.  In this paper we
recall the Gelfand integral and develop the calculus thereof,
including appropriate versions of the Fundamental Theorem of Calculus
and an integration by parts formula.  We then develop a related notion
of Gelfand weak (G-weak) derivative, and introduce the space
$W^{1,q}_{w*}(0,T;X^*)$ as the analogue of the usual space $W^{1,q}$
of weakly differentiable functions; in particular, our analysis shows
that functions from $W^{1,q}_{w*}(0,T;X^*)$ have absolutely continuous
representatives with weak derivatives in the space of weak* measurable
maps $L^q_{w*}(0,T; X^*)$, the dual of $L^{q'}(0,T;X)$.

Having developed the calculus of the Gelfand integral, we define the
notion of weak* solution.  Here weak* refers to the target space
rather than any type of convergence.  Suppose that the initial data
for \eqref{cl} is $u^0 \in X^*$.  We say that the function 
$u \in W^{1,q}_{w*}\big(0,T; X^*\big)$
is a \emph{weak* solution} to \eqref{cl}, if
\begin{equation}
\label{wssolintro}
  \bar{u}(t)-\bar{u}(s) = \wint_{s}^t \DD_x F(u(\tau))\, d\tau,
  \com{in} X^*,
\end{equation}
for $t$, $s \in [0,T]$, and $\bar u(0)=u^0$ in $X^*$.  Here $\bar{u}$
is the absolutely continuous representative of $u$ with values in
$X^*$, and we must make sense of the flux map $F: X^* \to X^*$ and
distributional derivative operator $\DD_x$, while implicitly requiring
$\DD_x  F(u(\cdot)) \in X^*$.  The formulation \eqref{wssolintro}
is equivalent to that of \eqref{BODE}, interpreting $u'$ as the
$G$-weak derivative of $u$.  Moreover, \eqref{wssolintro} provides a more
direct interpretation of the system \eqref{cl} than \eqref{weak},
essentially because the integration by parts is carried out implicitly
and abstractly in the definition of the spaces
$W^{1,q}_{w*}\big(0,T; X^*\big)$.

We note that, unlike the usual definition of weak solution, our
definition implicitly imposes a certain regularity on the solution
$u$, in that $t \to u(t)$ must be absolutely continuous in $X^*$.  Also,
although we require $\DD_xF(u(t))\in X^*$, this need only be weak*
measurable, rather than strongly measurable as a function of $t$. More
precisely, the mapping
\[
  t \to \DD_xF(u(t))  \com{belongs to}
  L^{q}_{w*}(0,T; X^*)=L^{q'}(0,T; X)^*.
\]
Note also that we allow any $1\le q\le\infty$, so we allow different
rates of growth of solutions.  When weak* solutions contain
discontinuties such as shocks, we are again required to rule out
inadmissible shocks via entropy considerations; when an entropy-flux
pair is available, the entropy inequality is again interpreted in the
weak* sense.

In principle, by changing the space $X$ of spatial test functions, we
could admit different classes of weak* solutions.  As an example,
classical $H^s$ solutions are obtained by setting $X=H^{-s}$, a large
space, while by restricting the class of test functions to a smaller
space, say $X=H^s$, we allow more solutions.  In \cite{MY2}, the
authors choose $X=C_0$, so that weak* solutions can include Radon
measures, which allows the treatment of vacuums in gas dynamics and
fractures in elasticity.  These weak* solutions are not weak solutions
as they include Dirac masses; on the other hand, the implicit
regularity of weak* solutions means that weak* solutions are more
restrictive than weak solutions.

In the present work, we primarily focus on weak* solutions
$u \in X^* \bigcap BV^n$, with $X=(C_0)^n$ and $X^*=M^n$, which we
call $BV$ weak* solutions.  Having defined this class of solutions, we
compare them to the more familiar entropy weak solutions.  We prove
that weak* solutions are distributional solutions, essentially because
the product test functions $\alpha(x)\,\beta(t)$ are dense in the set
$C^\infty_c(\RR^2)$.  In particular, a locally bounded weak* solution
is a weak solution.

Next, we prove that a weak solution with appropriate a.e.~bounds on
the spatial total variation is a weak* solution, and in particular,
any global solution obtained via Glimm's method or vanishing viscosity
is a global weak* solution with $q=\infty$.  This means that the
celebrated well-posedness theory for small $BV$ weak solutions
developed by Bressan's school applies without change to weak*
solutions.  Indeed, the $L^1$ stability of $BV$ weak* solutions
follows directly once existence of $BV$ weak* solutions is
established, see Theorem~\ref{w*wisweak}.

As an illustrative example of the use of $BV$ weak* solutions, we
derive the Rankine-Hugoniot jump conditions, by explicitly
differentiating jumps to get Dirac masses in the derivative, and
similarly express the entropy conditions.  We also show that the
quasilinear form of the equation is satisfied at any Lebesgue point of
a weak* solution.  Our derivation of the Rankine-Hugoniot conditions
is both easier and more general than the usual
derivation~\cite{Dafermos}.  As a further example, we express the
solution of the Riemann problem as a Bochner integral, which in this
case coincides with the Gelfand integral.

Finally, in an appendix, we consider the conditions that imply that an
abstract function be G-weak differentiable.  Indeed, we show that for
$1< q\le \infty$, the space $W^{1,q}_{w*}(0,T;X^*)$ is isometric to
the space of functions having bounded variation, introduced by Brezis
in~\cite{Brezis73}.

Our goal in this paper has been to introduce the calculus of the
Gelfand integral and to define the notion of weak* solution, showing
that it is consistent with the usual well-known class of $BV$ weak
solutions.  In the upcoming papers \cite{MY2,MY3,MY4}, we develop
these results, as follows.  In~\cite{MY2}, we extend the definition in
a natural way to rigorously allow the use of Dirac masses to deal with
vacuums and other discontinuities in the medium, and in~\cite{MY3} we
implement a front-tracking scheme for approximating solutions with
vacuum.  In \cite{MY4}, we extend these ideas to an abstract and
general framework in which we define weak* solutions for
multi-dimensional systems of balance laws.  It is hoped that the
availability of a wide variety of test function spaces $X$ will shed
light on the problems of instability and nonuniqueness of (strongly
measurable) weak solutions.

The paper is arranged as follows.  In Section 2, we recall the basic
facts we need and establish notation for the rest of the paper.  In
Section 3 we recall the Gelfand integral and develop the calculus
thereof: we do this for an abstract target space $X^*$, so we can use
these results unchanged in more general contexts.  In Section 4 we
use the Gelfand integral to define weak* solutions, and prove that
weak* solutions are distributional solutions, and weak solutions which
satisfy appropriate growth conditions are weak* solutions.  Finally, in
the appendix, we consider necessary and sufficient conditions for a
function to have a G-weak derivative.

\section{Preliminaries}

We begin by setting notation and recalling various facts about
functions, measures and derivatives that will be useful throughout the
paper.  For general references, we refer the reader to the books of
Royden, Cembranos \& Mendoza, Diestel \& Uhl, etc.~\cite{Royden2010,
  CM97, DU89}.

\subsection{Measures and Distributions}

Let $\Omega\subset\RR$ be an open set.  Denote the Borel
$\sigma$-algebra on $\Omega$ by $\mathcal{B}(\Omega)$, and Lebesgue
measure on the measure space $\big(\Omega, \mathcal{B}(\Omega)\big)$
by $\lambda$; we also write $dx=d\lambda=\lambda(dx)$.

For any distribution $T \in \mc D'(\Omega)$, we use the standard
pairing $\< T , \varphi \>$ to denote the action of $T$ on $\varphi
\in \mc D(\Omega)=C_c^{\infty}(\Omega)$.  Any integrable function $f\in
L^1_{loc}(\Omega)$ acts as distribution by
\[
  \< f, \varphi \> := \int f \, \varphi \; dx\,,
  \quad \varphi \in \mc D(\Omega).
\]
Let $M_{loc}(\Omega)$ denote the set of all sigma-finite Radon
measures on the measure space $(\Omega, \mc B(\Omega))$.  Any measure
$\nu \in M_{loc}(\Omega)$ is also a distribution via the action
\[
  \< \nu, \varphi \> := \int \varphi(x) \; \nu(dx)
   = \int \varphi(x) \; d\nu\,,\quad \varphi \in \mc D(\Omega)\,.
\]

We embed the integrable functions $L^1(\Omega)$ in the space of Radon
measures by means of the natural mapping $T:L^1(\Omega) \to M(\Omega)$,
$T(f)=:T_f$, defined by
\[
  T_f(E) = \int_{E} f(x) \; dx\,,\quad E\in \mc B(\Omega)\,.
\]
If $f\in L^1(\Omega)$, it is easy to see that $T_f$ is a Radon measure
which satisfies $\|T_f\|_{M(\Omega)}=\|f\|_{L^1(\Omega)}$, so that $T$
is norm preserving, and we can regard
\begin{equation}
  \label{embed}
  L^1_{loc}(\Omega) \subset M_{loc}(\Omega) \subset \mc D'(\Omega)\,,
  \quad
  L^1(\Omega) \subset M(\Omega) \subset \mc D'(\Omega)\,.
\end{equation}

Next, suppose that $F:\Omega\subset\RR\to\RR$, with $\Omega\subset\RR$
open.  We denote the classical derivative of $F$ at $x\in\Omega$ by
$\frac{dF}{dx}$, wherever this exists.  We adopt the convention that
if $\frac{dF}{dx}$ exists $\lambda$-almost surely, then we redefine
$\frac{dF}{dx}$ to be zero at all points $x$ at which the classical
derivative is not defined.  If $F\in L^1_{loc}(\Omega)$ has a weak
derivative, we denote the weak derivative by $F'(x)$ or $\D F\in
L^1_{loc}(\Omega)$.  Finally, we denote the distributional derivative
of $F$ by $\DD F \in \mc D' (\Omega)$.

\subsection{BV Functions}\label{BVfunc}

We recall a few basic facts about $BV$ functions.  The total variation
of a function $F:(a,b)\to\RR$ is
\[
  \Var(F, (a,b)) := \sup_{\mathcal{P}} \sum_{k=1}^N |F(x_{k+1})-F(x_k))|\,,
\]
where $\mathcal{P} := \{x_k\}_{k=1}^N$ ranges over the set of all
finite increasing sequences in $(a,b)$, and we say that $F$ is of
function of bounded variation on $(a,b)$ if $\Var(F,(a,b))<\infty$.
We denote the space of functions of bounded variation by
\[
   BV\big((a,b)\big) := \big\{ F:(a,b) \to R\,:\,\Var(F, (a,b))
   < \infty \big\}\,.
\]
We note that the usual convention is to treat $BV$ functions on a
closed interval $[a,b]$, but we find it convenient to consider an open
interval when studying their distributional derivatives.  Clearly, $F$
can be extended to the closed interval $[a,b]$ by continuity.  The
classical derivative $\frac{dF}{dx}$ exists Lebesgue almost surely,
the left and right limits $F(x\pm) := \lim_{z\to x\pm} F(z)$ exist for
all $x\in(a,b)$, and the set of all jump discontinuities
\[
  \mathcal{J}  = \{x \in (a,b): F(x+)-F(x-) \neq 0 \}
\]
is at most countable.  The right-continuous modification $F^*$ of $F$
satisfies $F^* = F_c + F_j + F_s$, where the functions
\begin{equation}\label{BVpart}
\begin{aligned}
F_{c}(x) &:= F(a) + \int_{a}^x \frac{dF}{dx}\; ds\,,\\
F_{j}(x) &:= \sum_{\bar{x} \in \mathcal{J}}
    H(x-\bar{x})\big(F(\bar{x}+)-F(\bar{x}-) \big)\,, \quad
     H(x) := \mc X_{x \geq 0}\,,\\
F_{s}(x) &:= F^*(x)-F_{c}(x)-F_{j}(x)\,,
\end{aligned}
\end{equation}
are the absolutely continuous, jump and singular continuous parts of
$F$, respectively.

Each $BV$ function gives rise to a finite signed Borel measure via
\begin{equation}
  \label{borel}
  \nu(y,x] := F(x+)-F(y+), \quad \forall x,\;y \in (a,b),
\end{equation}
and $\nu(a,b) := F(b)-F(a)$, with variation
$\|\nu\|_{M(a,b)} = |\nu|(a,b) = \Var(F,(a,b))$.
Applying \eqref{borel} to the decomposition \eqref{BVpart} of $F$
yields $\nu=\nu_c +\nu_{j}+\nu_{s}$, where $\nu_c$ is absolutely
continuous, $\nu_c\ll\lambda$, with $\lambda$-a.e.~Radon-Nikodym
derivative $\frac{d\nu_c}{d\lambda} = \frac{dF}{dx}$; $\nu_{j}$ is a
purely atomic singular measure $\nu_{j}\perp\lambda$; and
$\nu_{s}$ is a singular continuous measure, $\nu_s\perp\lambda$.
Moreover, we have the identities
\[
  \nu_c (A)  = \int_{A} \frac{dF_c}{dx} \; dx \com{and}
  \nu_{j} = \sum_{\bar{x} \in \mathcal{J}}
      \big( F(\bar{x}+)-F(\bar{x}-)\big)\,\delta_{\bar{x}}.
\]
These decompositions fully describe the distributional
derivative of the function $F\in BV(a,b)$, as follows.

\begin{lemma}
\label{DBV}
Using the above notation, $F \in BV(a,b)$ has distributional derivative
given by
\[
\begin{aligned}
  \DD F &= \DD F^* = \nu = \nu_c + \nu_j + \nu_s\,, \com{with} \\
  \DD F_c &= \nu_c\,, \quad \DD F_{j} = \nu_{j}\,, \com{and} \DD F_{s} = \nu_{s},
\end{aligned}
\]
so that
\[
\begin{aligned}
\< \DD F, \varphi \> &= \int \varphi\; \nu_c(dx)
   + \int \varphi \; \nu_{j}(dx) + \int \varphi \; \nu_{s}(dx)\\
& =  \int \varphi  \, \frac{dF}{dx} \; dx + \sum_{\bar{x} \in \mathcal{J}} \big( F(\bar{x}+)-F(\bar{x}-)\big)\,\varphi(\bar{x}) + \int \varphi \; \nu_{s}(dx)\,.
\end{aligned}
\]
\end{lemma}

\subsection{Products of Banach Spaces}

Let $X$ be a normed space, with norm denoted by $\|\cdot\|_X$.
Recall that the dual $X^*$ of $X$ has norm
\[
  \|\phi\|_{X^*} = \sup_{x \in X, x \neq 0} \frac{|\phi(x)|}{\|x\|_X}
   = \sup_{x \in X, x \neq 0} \frac{\<\phi,x\>}{\|x\|_X}\,,
\]
where $\<\phi,x\>$ denotes the pairing between $\phi$ and $x$, and
$X^*$ is a Banach space.  We recall the well-known property that this
pairing distinguishes elements in $X$ or $X^*$: that is, if
$\<\phi,x\> = 0$ for all $x\in X$, then $\phi=0$, and if $\< \phi, x
\> = 0$ for all $\phi\in X^*$, then $x=0$.

Because we are interested in systems, our solutions will consist of
vectors with values in $(X^*)^n$.  For any $X$, we equip the product
space $X^n$ with the ``Euclidean'' norm
\begin{equation}
  \label{norm}
  \|x\|_{X^n} := \Big(\sum_{i=1}^n \|x_i\|_X^2\Big)^{1/2}\,,
  \quad x = (x_1,\dots,x_n)\in X^n.
\end{equation}
We show that there is no ambiguity between $(X^*)^n$ and $(X^n)^*$.

\begin{lemma}
For any $n \geq 1$, the mapping $I: (X^*)^n \to (X^n)^*$, given by
\[
  [I(F)](x) = \sum_{i=1}^n F_i(x_i) \com{for}
  F=(F_1,\dots,F_i) \in (X^*)^n, \quad x\in X^n\,,
\]
is an isometric isomorphism, $\|I(F)\|_{(X^n)^*}=\|F\|_{(X^*)^n}$.
\end{lemma}

\begin{proof}
  The map $I$ is trivially injective and onto; to show it is an
  isometry, first note that
\[
   \Big|[I(F)](x)\Big| = \Big|\sum_{i=1}^nF_i(x_i)\Big|
   \leq \sum_{i=1}^n \|F_i\|_{X^*}\,\|x_i\|_X \le \|F\|_{(X^*)^n}\|x\|_{X^n},
\]
by Cauchy-Schwarz in $\RR^n$, so that
$\|I(F)\|_{(X^n)^*} \leq \|F\|_{(X^*)^n}$.  Next, given any
$F \in (X^*)^n$ and $\eps >0$, we can find $\bar{x}_i \in X$,
$i=1,\dots,n$ such that each $\|\bar{x}_i\|=1$ and
\[
  \Big(\sum_{1=1}^n\Big|F_i(\bar{x}_i) - \|F_i\|_{X^*}\Big|^2 \Big)^{1/2} \leq
  \sqrt n\,\eps.
\]
Now set
\[
  \bar{y} = \|F\|_{(X^*)^n}^{-1}\,\Big( \,\bar{x}_1\|F_1\|_{X^*}\,,\,
            \dots\,,\,\bar{x}_n\|F_n\|_{X^*} \Big),
\]
so that $\|\bar{y}\|_{X^n} = 1$, and
\[
   [I(F)](\bar{y}) - \|F\|_{(X^*)^n} = \|F\|_{(X^*)^n}^{-1}
   \sum_{i=1}^n \|F_i\|_{X^*}\big(F_i( \bar{x}_i )-\|F_i\|_{X^*}\big)
   \ge - \sqrt n\,\eps \,,
\]
again using Cauchy-Schwarz.  Since $\eps$ is arbitrary, the result follows.
\end{proof}

\begin{corollary}\label{radreprn}
  Let $\Omega \subset \RR^d$ be an open set, and $C_0(\Omega)$ denote
  the closure of $C_c(\Omega)$ under the sup-norm.  For $n \geq 1$,
  the mapping $T: M(\Omega)^n \to (C_0(\Omega)^n)^*$, given by
\[
  \< T(\mu) , \varphi \> = \int_{\Omega} \varphi(x) \cdot \mu(dx)
   :=  \sum_{i=1}^n \int \varphi_i(x) \; \mu_i(dx)\,,
\]
is an isometric isomorphism,
\[
  \|T(\mu)\|_{(C_0(\Omega)^n)^*}
  = \|\mu\|_{M^n(\Omega)}
  = \Big(\sum_{i=1}^n(|\mu_i|(\Omega))^2\Big)^{1/2}.
\]
\end{corollary}

\begin{proof}
  Since $\Omega \subset \RR^d$ is locally compact, the proof for $n=1$
  follows from the Riesz-Markov and Riesz Representation
  theorems~\cite{Royden2010}; the general case follows directly.
\end{proof}

\subsection{Measurability and Equivalence of Functions}

Because we want to integrate in Banach spaces, we recall several
notions of measurability and note the relationships between them.  We
follow the conventions used in the book of Diestel and Uhl
\cite{DU89}.  Unless explicitly stated otherwise, $f$ is a map of a
single variable with values in a Banach space, $f:[0,T]\to X$, and the
Borel $\sigma$-algebra is always assumed.

The map $f$ is \emph{measurable}, or \emph{$\BB$-measurable}, if
$f^{-1}(A)\in \BB := \mc B([0,T])$ for each $A \in \mc B(X)$\,.

The map $f(t)$ is called \emph{$\lambda$-essentially separably valued}
if there exists $N\in \BB $ with $\lambda(N)=0$ and a countable set $H
\subset X$ such that $f([0,T] \backslash N) \subset \overline{H}$; $f$
is \emph{separably valued} if $f([0,T]) \subset \overline{H}$.

The map $f$ is \emph{simple} if there exist vectors
$u_1,\ u_2,\dots,u_n \in X$ and sets $E_1,\ E_2,\dots,E_n \in \BB$ such
that
\[
f(t)=\sum_{i=1}^n u_i \; \mc X_{E_i}(t)
\]
where $\mc X_{E_i}$ is the indicator function of $E_i$.

The map $f$ is \emph{$\lambda$-measurable}, or \emph{strongly
  measurable}, if there is a sequence of simple functions $f_n:
[0,T]\to X$ such
that
\[
  \lim_{n\to \infty}\| f_n(t) - f(t) \| = 0
  \quad \mbox{$\lambda$-almost surely on  $[0,T]$.}
\]

It is not hard to show that a map is $\lambda$-measurable if and only
if it is both $\BB$-measurable and $\lambda$-essentially separably
valued.  In particular, if the space $X$ is separable, then
$\lambda$-measurability and $\BB$-measurability are equivalent.

The map $f:[0,T]\to X$ is \emph{norm-measurable} if its norm
$\| f(t) \|:[0,T]\to \RR$ is $\lambda$-measurable (i.e.~Lebesgue
measurable).

If $f:[0,T]\to X$ is $\lambda$-measurable, it is norm-measurable.
Moreover, if $f_n:[0,T]\to X$ is a sequence of $\lambda$-measurable
functions which almost surely converge in norm, $\| f_n(t) - f(t) \|
\to 0$ as $n\to\infty$ $\lambda$-almost surely on $[0,T]$, then the
limit $f$ is $\lambda$-measurable.

We now define a map $f: [0,T]\to X$ to be \emph{weakly
  measurable}, or \emph{scalar measurable}, if for each $\phi
\in X^*$, the function
\[
  \phi(f(\cdot)) = \< \phi, f(\cdot) \>: [0,T]\to \RR
\]
is $\lambda$-measurable, or equivalently $\BB$-measurable.

A measurable map is weakly measurable, and the converse holds if $f$
is $\lambda$-essentially separably valued; this is the statement of
Pettis' Measurability Theorem; see~\cite{DU89}.  It follows
immediately that if $X$ is separable, then $f$ is $\lambda$-measurable
if and only if it is weakly measurable.  Pettis's theorem also
yields the fact that $f: [0,T]\to X$ is $\lambda$-measurable if and
only if it is the $\lambda$-almost everywhere uniform limit of a
sequence of countably valued $\lambda$-measurable functions.

The map $f:[0,T]\to X^*$ is \emph{weak*-measurable} if for
each $x \in X$, the function
\[
  [f(\cdot)](x)=\< f(\cdot), x \>: [0,T] \to \RR
\]
is $\lambda$-measurable.

In view of the standard embedding $X \subset X^{**}$, it follows that
an $X^*$-valued weak measurable map $f: [0,T]\to X^*$ is
weak* measurable.

We now introduce notions of equivalence of functions.  Suppose we are
given $f$, $g:[0,T]\to X$ and $f^*$,~$g^*:[0,T]\to X^*$.

The maps $f$ and $g$ are \emph{$\lambda$-equivalent} if $f(t)=g(t)$
$\lambda$-a.e.

The maps $f$ and $g$ are \emph{weakly equivalent} if, for
all $\phi\in X^*$, and for $\lambda$-a.e.~$t$, we have 
$\< \phi,f(t) \> = \< \phi, g(t) \>$.

The maps $f^*$ and $g^*$ are \emph{weak* equivalent} if, for
all $x\in X$, and for $\lambda$-a.e.~$t$, we have 
$\< f^*(t),x \> = \< g^*(t),x \>$.

It is clear that $\lambda$-equivalent functions are weakly equivalent,
and that weakly equivalent functions with values in $X^*$ are weak*
equivalent.  In fact, the converse of this is true, provided the
functions are strongly measurable; see~\cite{DU89}.

\begin{lemma}
  Suppose that $f$, $g:[0,T]\to X$ and $f^*$, $g^*:[0,T]\to X^*$ are
  $\lambda$-measurable.  If $f$ and $g$ are weakly equivalent, then
  they are $\lambda$-equivalent, and if $f^*$ and $g^*$ are weak*
  equivalent, then they are $\lambda$-equivalent.
\end{lemma}

In this lemma, the requirement that the functions be strongly
measurable is essential.  We will later consider non-measurable
functions that are weak or weak* measurable, in which cases the
corresponding weaker notions of equivalence are appropriate.

\subsection{The Bochner Integral}

The integral of a simple function, given by
$h(t)=\sum_{i=1}^N u_i\, \mc X_{E_i}: [0,T] \to X$, is defined in the
obvious way,
\[
  \int_0^T h(t) \; \lambda(dt) = \int_0^T h(t) \; dt
   := \sum_{i=1}^N u_i \; \lambda(E_i)\,\in X.
\]
We say that a $\lambda$-measurable function $ f: [0,T] \to X$ is
\emph{Bochner integrable}, or \emph{summable}, if there exists a
sequence of simple functions $\{h_n\}_{n \geq 1}$ such that the
Lebesgue integral $\int_0^T \| h_n - f\| \; dt\to 0$ as $n\to\infty$.
It follows that if $f$ is summable, its integral over any $E \in \BB$
exists in the space $X$,
\[
  \int_E f(t) \; dt := \lim_{n \to \infty} \int_E h_n(t) \; dt \,.
\]

The calculus of the Bochner integral is well known, and many of the
usual theorems of the Lebesgue integral translate directly to
statements on the Bochner integral.  These include a Bochner dominated
convergence theorem and a Bochner-Lebesgue differentiation theorem,
which states that for $\lambda$-almost all $s \in [0,T]$,
\begin{equation}
\label{BLDT}
\begin{aligned}
  \lim_{h \to 0} \frac{1}{h} \int_s^{s+h} \| f(t) - f(s) \|\; dt &= 0
\com{and}  \\
  \lim_{h \to 0} \frac{1}{h} \int_s^{s+h} f(t) \; dt &= f(s)\,.
\end{aligned}
\end{equation}

Hille's theorem states that Bochner integration commutes with closed
linear operators, and in particular implies that if $f:[0,T]\to X$
and $f^*:[0,T] \to X^*$ are Bochner integrable, then for each $E \in
\BB$, $x\in X$ and $\phi\in X^*$,
\[
\begin{aligned}
  \Big\langle \phi, \int_E f(t)\, dt \Big\rangle &=
     \int_E \< \phi, f(t) \> \; dt \com{and} \\
  \Big\langle \int_E f^*(t)\, dt, x \Big\rangle &=
     \int_E \< f^*(t), x \> \; dt \,.
\end{aligned}
\]

The Bochner integral allows us to introduce $L^p$ spaces as a natural
generalization of the usual space $L^p(0,T)$ of real valued
functions.  For each $1 \leq p \leq \infty$, the Banach space
$L^p(0,T; X)$ consists of all $\lambda$-equivalence classes of Bochner
integrable functions $f:[0,T]\to X$ such that
$\|f\|_{L^p(0,T;X)}<\infty$, where
\[
  \| f \|_{L^p(0,T;X)} :=
   \begin{cases}
    \ \left( \int_0^T \|f(t)\|^p\; dt\right)^{1/p}\,,
        & 1\le p < \infty\,,\\[2pt]
    \ \esssup_{t\in[0,T]} \| f(t)\|\,, & p = \infty\,.
  \end{cases}
\]

Now let $X,Y$ be Banach spaces with $Z=X \bigcap Y$ nonempty, and
suppose $f\in L^1(0,T;X)$.  We say that $g \in L^1(0,T; Y)$ is the
\emph{Bochner weak derivative} or B-weak derivative of $f$, written
$f'(t)=g(t)$, provided that $f$, $g$ are Bochner integrable as
functions from $[0,T] \to Z$ and
\[
  \int_0^T f(t) \, \varphi'(t)\; dt
   = - \int_0^T g(t) \, \varphi(t) \; dt
\]
for all scalar functions $\varphi \in C^{\infty}_c(0,T)$.  For
$X^*$-valued functions, B-weak derivatives are fully determined by
their actions as functionals: that is, for $f$, $g\in L^1(0,T;X^*)$, ${f}'=g$
if and only if for every $x \in X$ and
$\varphi \in C^{\infty}_c(0,T)$,
\[
  \int_0^T \< f(t),x \> \varphi'(t) \, dt
    = - \int_0^T \< {g}(t) , x\> \varphi(t) \, dt.
\]

For $1 \leq p \leq \infty$, the Sobolev space $W^{1,p}\big(0,T;X\big)$
is the set of all functions $f\in L^p\big(0,T;X\big)$ with
$f'\in L^p\big(0,T;X)$, with norm
\[
  \|f\|_{W^{1,p}(0,T;X)} :=
\begin{cases}
\ \big(\int_0^T \|f(t)\|^p+\|f'(t)\|^p \;dt\big)^{1/p} \,,
    & 1\leq p < \infty,\\[2pt]
\ \esssup_{t\in[0,T]}\big( \|f(t)\| + \|f'(t)\| \big)\,, & p=\infty\,.
\end{cases}
\]

For our purposes, the most important feature of functions in
$W^{1,p}(0,T;X)$ is the Fundamental Theorem of Calculus, which follows
from the fact that if $f\in L^{p}(0,T;X)$, the integral
$\int_0^tf(s)\,ds$ is absolutely continuous as a function of $t$.

\begin{theorem}\label{BFTC}
  Let $f\in W^{1,p}(0,T;X)$ for some $1\leq p \leq \infty$, and define
  $\bar f(t) = f(0) + \int_0^t f'(s)\, ds$.  Then $\bar f \in
  C\big([0,T]; X \big)$ is a.e.~differentiable and $\bar f = f$ a.e..
  The function $g:[0,T] \to X$ given by
\[
g(s) := \begin{cases}
  \ \lim_{h\to 0}\frac{1}{h}\big(\bar{f}(s+h)-\bar{f}(s)\big)\,,
   & \text{if the limit exists,}  \\
  \ 0\,,& \text{otherwise,}
\end{cases}
\]
is strongly measurable and satisfies $g(t) = f'(t)$ $\lambda$-a.e.~on
$[0,T]$.  If $p=\infty$, then $\bar{f}$ is Lipschitz continuous with
\[
  \big\|\bar{f}(t)-\bar{f}(s)\big\| \leq L\, |t-s|  \quad \mbox{with} \quad
  L := \esssup_{t\in[0,T]}\|f'(t)\|\,.
\]
\end{theorem}

The following form of H\"older's inequality holds: if $f\in
L^p(0,T;X)$ and $g^* \in L^q(0,T;X^*)$, $p^{-1}+q^{-1}=1$, then
$\< g^*, f \> :[0,T] \to \RR$ is integrable and satisfies
\[
  \int_0^T \big| \< g^*(s), f(s) \> \big|\; ds \leq
  \big\| f \big\|_{L^p(0,T;X)}\; \big\| g^* \big\|_{L^q(0,T;X^*)}\,.
\]
As a consequence, we can integrate by parts: given $f\in W^{1,p}(0,T;X)$
and $g^* \in W^{1,q}(0,T;X^*)$, we have $\<g^*,f\> \in W^{1,1}(0,T)$,
and
\begin{equation}\label{intbyparts}
  \int_0^T \< {g^*}'(t), f(t) \> \; dt
  = \<\bar{g}^*(s),\bar{f}(s)\>\, \Big|_0^T
  - \int_0^T \< {g^*}(t), f'(t) \>\;dt
\end{equation}
where $\bar{f}$, $\bar{g}$ are the continuous representatives of $f$
and $g$, respectively.

\section{The Gelfand Integral}

It is sometimes too restrictive to consider only the Bochner integral
of strongly measurable functions with values in a Banach space.  Here
we recall the Gelfand integral, which is better suited to our
purposes.  Throughout this section we fix a Banach space $X$ and let
$X^*$ denote its dual.  The Gelfand integral is similar to the Dunford
integral, but takes its values in the space $X^*$.  For a detailed
discussion see Tulcea and Tulcea~\cite{IT1969}, Diestel and
Uhl~\cite{DU89}, or Cembranos and Mendoza \cite{CM97}.

Suppose that we are given a weak*-measurable function $\Phi:[0,T]\to
X^*$, and suppose also that
\[
  \< \Phi(\cdot), x \> \in L^1(0,T) \com{for all} x \in X\,.
\]
For a given $E \in \BB = \mc B([0,T])$, we define the map
$T_E:X\to L^1(0,T)$ by
\[
  T_E(x) = \< \Phi(\cdot), x \>\,\chi_E(\cdot) \in L^1(0,T)\,.
\]
It is clear that $T_E$ is linear, and if $x_n\to x$ and $T_E(x_n)\to
y$ in $L^1$, then by the Riesz-Fischer theorem, a subsequence
$T_E(x_{n_k})(s)\to y(s)$ a.e., while also $T_E(x_n)(s) \to \<
\Phi(s), x \>\,\chi_E(s)$ for all $s\in[0,T]$.  It follows that $y\in
L^1(0,T)$, so $T_E$ is closed, and further, by the closed graph
theorem, it is bounded, so we can write $\| T_E(x) \|_{L^1} \le
\|T_E\|\,\|x\|$ for all $x\in X$.  Since integration is a bounded
linear operator of $L^1$ into $\RR$, it follows that the map
\[
  x \mapsto \int_0^T T_E(x)(s)\;ds = \int_E\< \Phi(s), x \>\;ds
\]
is a bounded linear functional on $X$, so defines an element of the
dual $X^*$.  

\begin{definition} Let $\Phi:[0,T] \to X^*$ be a weak* measurable
  function such that $\< \Phi(\cdot), x \> \in L^1(0,T)$ for every
  $x \in X$.  The Gelfand integral of $\Phi$ over measurable set
  $E \subset (0,T)$, denoted by $\wint_E \Phi(s)\; ds$, is an element
  of $X^*$ defined by
\begin{equation}
  \label{G-int}
  \Big\langle \wint_E \Phi(s)\;ds , x \Big\rangle =
     \int_E \< \Phi(s),x\>\;ds \com{for all} x\in X\,.
\end{equation}
\end{definition}

\subsection{The Spaces $L^q_{w*}(0,T;X^*)$}

It is clear from \eqref{G-int} that if two weak*-measurable functions
are weak*-equivalent, then their Gelfand integrals coincide.  We are
thus led to consider equivalence classes,
\[
[\Psi] = \big\{\widetilde{\Psi}:[0,T]\to X^*:
   \mbox{$\Psi$ and $\widetilde{\Psi}$ are weakly*-equivalent} \big\}\,.
\]
We now wish to describe the sets of (equivalence classes of)
weak*-measurable functions which are $L^q$ Gelfand integrable, and to
give an appropriate norm.  Note that even if $\Phi:[0,T]\to X^*$
is weak*-measurable, in general its norm $\|\Phi\|:[0,T]\to\RR$ is not
measurable.  However, the theory of liftings provides the following
lemma,~\cite{CM97,IT1969}.

\begin{lemma}
  \label{normmeas}
  Each Gelfand integrable function $\Psi:[0,T]\to X^*$ is
  weak*-equivalent to a map $\hat\Psi:[0,T]\to X^*$ whose norm
  $\|\hat\Psi\|:[0,T]\to \RR$ is measurable.  Moreover, if $\Psi_1$ and
  $\Psi_2$ are weak*-equivalent and norm-measurable, then
  $\|\Psi_1\|=\|\Psi_2\|$ almost everywhere.
\end{lemma}

It is easy to show that if $X$ is separable, all weak*-measurable
functions are norm-measurable.

We are now in a position to describe the $X^*$ valued Gelfand $L^q$
spaces, for $1\le q\le\infty$.  Given an equivalence class $[\Psi]$ of
Gelfand integrable functions, set
\[
  \qn{[\Psi]}_q := \inf\big\{ \|g\|_{L^q(0,T)}\  :\
  \|{\hat\Psi}(t)\|\leq g(t)\ \text{$\lambda$-a.e.} \big\}\,,
\]
where $\hat\Psi\in[\Psi]$ is a norm-measurable element of the
equivalence class.  It follows that $\qn{\cdot}_q$ is a norm, and we let
$L^{q}_{w*}(0,T;X^*)$ be the space of equivalence classes $[\Psi]$ of
finite norm,
\[
  L^{q}_{w*}(0,T;X^*) := \big\{ [\Psi] :
     \qn{[\Psi]}_q < \infty \big\}\,.
\]

It is not difficult to show that $L^q_{w*}(0,T;X^*)$ is a Banach space
and that the trivial inclusion
\[
L^q(0,T;X^*) \subset L^q_{w*}(0,T;X^*) \com{via} f \mapsto [f],
\]
is a norm-preserving isomorphism.  In particular, a Bochner integrable
map $f:[0,T]\to X^*$ is also Gelfand integrable, and when they both
exist, the integrals coincide.  Also, if
$\hat\Psi\in[\Psi]\in L^q_{w*}(0,T;X^*)$ is norm-measurable, then
$\|\hat\Psi\|\in L^q(0,T)$ and
\[
  \qn{[\Psi] }_q
  = \|\;\|\hat\Psi(\cdot)\|\;\|_{L^q(0,T)}\,.
\]
Similarly,
\[
  \int_0^T \big|\< \hat\Psi(t), x \>\big|^q \; dt
  \leq \int_0^T \|\hat\Psi(t)\|^q \,\|x\|^q \; dt
  = \qn{[\Psi]}_q^q\; \|x\|^q,
\]
and we conclude that for any $x\in X$ and
$\wt\Psi\in[\Psi] \in L^q_{w*}(0,T;X^*)$, we have
\begin{equation}
  \label{w*integrbnd}
  \< \wt{\Psi}(\cdot),x \> \in L^q(0,T) \com{with}
  \big\| \< \wt{\Psi}(\cdot) , x \> \big\|_{L^q(0,T)} \leq
  \qn{[\Psi]}_q \|x\|\,.
\end{equation}

It turns out that $L^q_{w*}(0,T;X^*)$ is the dual space of
$L^p(0,T;X)$, provided $1 \leq p < \infty$ and $p^{-1}+q^{-1}=1$; see
\cite{CM97} for details.

\begin{theorem}\label{dualLpX}
For each $\widetilde\Psi\in [\Psi] \in L^q_{w*}(0,T;X^*)$ and
$f \in L^p(0,T;X)$, the function
\[
  \< \widetilde{\Psi}(\cdot),f(\cdot) \> \in L^1(0,T)\,,
\]
and is independent of representative $\widetilde\Psi$.  Moreover, the
natural linear map $[\Psi]\to I([\Psi])\in L^p(0,T;X)^*$, given by
\[
  \< I[\Psi], f \> := \int_0^T \< \widetilde{\Psi}(s),f(s) \>\;ds\,,
   \quad f\in L^p(0,T;X),
\]
is an isometric isomorphism of $L^q_{w*}(0,T;X^*)$ onto $L^p(0,T;X)^*$.
\end{theorem}

\subsection{Absolute Continuity}

Although the spaces $L^q_{w*}(0,T;X^*)$ can be regarded as well known,
to the authors' knowledge the calculus of Gelfand integrable functions
has not previously been developed.  Here we develop this calculus in
parallel with the well-known calculus of the Bochner integral.  We
begin by proving absolute continuity of the Gelfand integral.

\begin{theorem}\label{w*intprop}
  Fix $1\le q\le\infty$, and suppose $[\Phi] \in L^q_{w*}(0,T;X^*)$
  with $\wt{\Phi} \in [\Phi]$.  The linear map $\Psi(t):[0,T] \to X^*$
  defined by
\[
\begin{aligned}
  \Psi(t) &= \wint_0^t \wt{\Phi}(s) \; ds\,, \com{that is}\\
  \< \Psi(t),x \> &= \int_0^t \< \wt\Phi(s), x \>\;ds\,,
     \quad x \in X\,,
\end{aligned}
\]
is absolutely continuous.  Also, the total variation function
$\Var_{\Psi}:[0,T]\to\RR$, defined by
\[
  \Var_{\Psi}(t) := \Var\big(\Psi,[0,t]\big)
   = \sup_{\mc P_{[0,t]}}\sum \| \Psi(t_n) - \Psi(t_{n-1}) \|\,,
\]
the supremum being over partitions $\mc P_{[0,t]} = \{0=t_0<t_1\cdots
<t_n=t\}$, is absolutely continuous, with pointwise derivative
$\frac{d \Var_{\Psi}}{dt}\in L^q(0,T)$, and we have
\begin{equation}\label{w*var}
   \|\Psi(t)-\Psi(s)\| \leq  \Var_{\Psi}(t) - \Var_{\Psi}(s)
   = \int_s^t \frac{d\Var_{\Psi}}{dt}(\tau)\; d\tau
   \leq \int_s^t \| \hat\Phi (\tau) \| \; d\tau
\end{equation}
for all $0\leq s \leq t \leq T$, for $\hat\Phi\in[\Phi]$
norm-measurable.
\end{theorem}

\begin{proof}
  Suppose that $\hat\Phi\in[\Phi]$ is norm-measurable, and let
  $E\subset[0,T]$ be a union of disjoint intervals,
  $E=\bigcup_{n=1}^N(a_n,b_n)$.  Then $\|\hat\Phi\| \in L^1(0,T)$, and
  we have
\[
\begin{aligned}
  \sum_{n=1}^N \|\Psi(b_n)-\Psi(a_n)\|
  &= \sum_{n=1}^N \sup_{\|x\|=1}\< \Psi(b_n)-\Psi(a_n), x \>\\
  &= \sum_{n=1}^N \sup_{\|x\|=1}
         \int_{a_n}^{b_n}\< \hat\Phi(s), x \>\; ds \\
& \leq \sum_{n=1}^N \int_{a_n}^{b_n}\|\hat\Phi(s)\| \; ds
   = \int_E\|\hat\Phi(s)\| \; ds\,,
\end{aligned}
\]
and absolute continuity follows from that of the Lebesgue integral.
Similarly, for any $0\leq s \leq t \leq T$,
\[
\Var_{\Psi}(t)-\Var_\Psi(s)  \leq \int_s^t \|\hat\Phi(\tau)\| \; d\tau\,,
\]
so $\Var_{\Psi}(t)$ is absolutely continuous.  Thus the derivative
$\frac{d\Var_{\Psi}}{dt}$ is defined $\lambda$-a.e., and
\[
  \Var_{\Psi}(t)-\Var_\Psi(s)
  = \int_s^t \frac{d\Var_{\Psi}}{dt}(\tau) \, d\tau
  \quad \mbox{for all} \quad 0\leq s \leq t \leq T\,.
\]
It follows that $0 \leq \frac{d \Var_{\Psi}}{dt} \leq \|\hat\Phi(t)\|$
$\lambda$-a.e.~on $[0,T]$, and thus $\frac{d \Var_{\Psi}}{dt} \in
L^q(0,T)$.
\end{proof}

The Lebesgue Differentiation Theorem implies that for all $x\in X$,
\[
  \lim_{h \to 0} \frac{1}{h} \int_t^{t+h}
     \< \wt{\Phi}(s) , x \>\; ds
  = \< \wt{\Phi}(t), x \>\,   \quad \mbox{ $\lambda$-a.e.}\,.
\]
Our next task is to extend this to get a Gelfand-Lebesgue
Differentiation Theorem.  Fix $1 \leq q \leq \infty$ and
$p^{-1}+q^{-1}=1$.

\begin{theorem}\label{w*LDT}
Suppose that $[\Phi] \in L^{q}_{w*}(0,T;X^*)$ and $\wt\Phi\in[\Phi]$.
Then for each $f \in L^{p}(0,T;X)$,
\begin{align}\label{w*LDT2}
\lim_{h \to 0} \frac{1}{h} \int_t^{t+h} \<\wt{\Phi}(s) , f(t) \>\; ds
   &= \<\wt{\Phi}(t), f(t) \> \quad \mbox{ $\lambda$-a.e.}\,,\com{and}\\[2pt]
\label{w*LDT3}
\lim_{h \to 0} \frac{1}{h} \int_t^{t+h} \<\wt{\Phi}(s) , f(s) \>\; ds
   &= \<\wt{\Phi}(t), f(t) \> \quad \mbox{ $\lambda$-a.e.}\,.
\end{align}
If $X$ is separable then for almost all $t \in (0,T)$
\[
 \frac{1}{h} \wint_t^{t+h} \wt{\Phi} (s) \; ds
  \;\; \overset{w*}\longrightarrow \;\; \wt{\Phi} (t)
\quad \text{as \;\; $h \to 0$.}\\
\]
\end{theorem}

\begin{proof}
  Since $f \in L^p(0,T;X)$, it is strongly measurable,
  and without loss of generality we may assume that it is separably
  valued.  Let $\{x_n\}$ be a countable dense subset of $f([0,T])$ and
  let $\hat\Phi \in [\Phi]$ be norm-measurable.  Define the set
\[
\begin{aligned}
  E := &\bigg\{t:\,\frac{1}{h}\int_t^{t+h}\|\hat\Phi(s)\|\;ds
          \to \|\hat\Phi(t)\|\bigg\}\\
       &\quad\bigcap \mathop{\bigcap}_{n=1}^\infty
       \bigg\{t:\,\frac{1}{h}\int_t^{t+h}\<\wt{\Phi}(s),x_n\>\;ds
          \to \< \wt{\Phi}(t), x_n \>\,\bigg\}
\end{aligned}
\]
where the limits are taken as $h\to0$.  By the Lebesgue
Differentiation Theorem, $E$ has full measure,
$\lambda\big([0,T]\setminus E\big)=0$.  For $t\in E$ and any $n$, we have

\[
\begin{aligned}
 \Big|\frac{1}{h} \int_t^{t+h} &\< \wt{\Phi}(s) , f(t) \> \; ds
      - \< \wt{\Phi}(t), f(t)\> \Big|
\\& \leq \Big|\frac{1}{h} \int_t^{t+h}
       \< \wt{\Phi}(s) , f(t)-x_n \> \;ds \Big|\\
&  \qquad
  + \Big|\frac{1}{h} \int_t^{t+h} \< \wt{\Phi}(s) , x_n \> \; ds
          - \< \wt{\Phi}(t), x_n\> \Big|
  + \Big| \< \wt{\Phi}(t), x_n-f(t)\> \Big|  \\
& \leq \Big|\frac{1}{h} \int_t^{t+h} \|\hat\Phi(s)\|\; ds
           \Big|\;\|f(t)-x_n\|  \\
& \qquad + \Big|\frac{1}{h} \int_t^{t+h} \< \wt{\Phi}(s) , x_n \> \; ds - \< \wt{\Phi}(t), x_n\> \Big| + \|\wt{\Phi}(t)\| \; \|x_n-f(t)\|\,.
\end{aligned}
\]
Since $t \in E$ and $\{x_n\}$ is dense in $f([0,T])$, the right hand
side can be made arbitrarily small, and \eqref{w*LDT2} follows.

Now set $B_n = \{ t: \|\hat\Phi(t)\|<n \}$.  Clearly $B_n$ is nonempty
for all $n$ beyond some $N \geq 1$, and $\bigcup_nB_n$ has full measure.
For any $n \geq N$ and $t \in E \bigcap B_n$, using \eqref{BLDT},
we get
\[
\begin{aligned}
\limsup_{h \to 0}\Big|\frac{1}{h}
   \int_t^{t+h} & \< \wt{\Phi}(s) , f(s) - f(t) \> \; ds \Big|
\\& \leq \limsup_{h \to 0}\Big|\frac{1}{h} \int_t^{t+h}
        \|\hat\Phi(s)\|\; \| f(s)-f(t)\| \; ds \Big|  \\
& \leq n \;\limsup_{h \to 0}
   \Big|\frac{1}{h} \int_t^{t+h} \| f(s)-f(t)\| \; ds \Big| = 0,
\end{aligned}
\]
and thus \eqref{w*LDT3} follows from \eqref{w*LDT2}.

Now suppose that $X$ is separable and let $\{x_n\}_{n \geq 1}$ be dense
in $X$.  Define $E$ as above and fix $t\in E$.  Then there is some
$\delta>0$ such that for all $0<|h|<\delta$,
\[
  \frac{1}{h} \int_{t}^{t+h} \|\hat\Phi(s)\|  \, ds \leq \|\hat\Phi(t)\|+1\,.
\]
For each such $h$, define $\phi_{t}(h)\in X^*$ by
\[
  \phi_{t}(h) := \frac{1}{h}\wint_t^{t+h}\wt{\Phi}(s)\;ds,
\]
so that, for every $x \in X$,
\[
  \big|\<\phi_t(h), x\>\big| \leq
   \frac{1}{h} \int_t^{t+h} \|\hat\Phi(s)\|  \; ds \;\|x\|
 \le  \big(\|\hat\Phi(t)\| +1\big)\; \|x\|,
\]
which yields $\|\phi_t(h)\|_{X^*} \leq \|\hat\Phi(t)\| +1$ for each
$0<|h|<\delta$.  Now for any fixed $x\in X$ and each $n$, we have
\[
\begin{aligned}
  \big|\<\phi_t(h)-\wt{\Phi}(t),x\>\big| &\leq
   \big|\<\phi_t(h)-\wt{\Phi}(t),x_n\>\big|
   + \big|\<\phi_t(h),x-x_n\>\big|
   + \big|\<\wt{\Phi}(t),x_n-x\>\big|  \\
& \leq \big|\<\phi_t(h)-\wt{\Phi}(t),x_n\>\big|
   + \big(\|\hat\Phi(t)\|+1+\|\wt{\Phi}(t)\|\big)\; \|x-x_n\|\,,
\end{aligned}
\]
and since $t\in E$, for every $n \geq 1$ we have
\[
  \limsup_{h \to 0} \big|\<\phi_t(h)-\wt{\Phi}(t),x\>\big|
    \leq \big(\|\hat\Phi(t)\|+1+\|\wt{\Phi}(t)\|\big)\;\|x-x_n\|.
\]
Since $\{x_n\}$ is dense, we have
$\<\phi_t(h),x\>\to\<\wt{\Phi}(t),x\>$, and since $x$ is arbitrary,
the result follows.
\end{proof}

\subsection{The Spaces $W^{1,q}_{w*}(0,T;X^*)$}

Suppose that $\Psi$, $\Phi:[0,T]\to X^*$ are weak* measurable, and
$[\Psi]$, $[\Phi]\in L^1_{w*}(0,T;X^*)$.  We say that $\Phi$ is the
\emph{Gelfand weak derivative} or G-weak derivative of $\Psi$, written
$\Psi'(t)=\Phi(t)$ or $[\Phi]=[\Psi']$, if
\begin{equation}
\begin{aligned}
  \wint_0^T \Psi(t) \, \varphi'(t)\; dt
   &= - \wint_0^T \Phi(t)\, \varphi(t)\; dt, \com{that is}\\
  \int_0^T \< \Psi(t),x\> \, \varphi'(t)\; dt
   &= - \int_0^T \< \Phi(t),x\> \, \varphi(t) \; dt
   \com{for all} x\in X,
\end{aligned}
\label{w*wd}
\end{equation}
for all scalar functions $\varphi \in C^{\infty}_c(0,T)$.

We now define the space $W^{1,q}_{w*}\big(0,T;X^*\big)$, for
$1 \leq q \leq \infty$, to be the set of weak* equivalence classes
$[\Psi]\in L^q_{w^*}\big(0,T;X^*\big)$ having G-weak derivative
$[\Psi'] \in L^q_{w*}\big(0,T;X)$, with norm
\[
  \qn{ [\Psi] }_{W^{1,p}_{w*}(0,T;X)} :=
  \begin{cases}
    \big(\int_0^T (\|\hat\Psi(t)\|^q+\|\hat\Psi'(t)\|^q)\;dt \big)^{1/q} \,,
      & 1\leq q < \infty\\
    \esssup_{t\in[0,T]}\big( \|\hat\Psi(t)\| + \|\hat\Psi'(t)\| \big)\,,
      &  q=\infty\,,
  \end{cases}
\]
for $\hat\Psi \in [\Psi]$, $\hat\Psi' \in [\Psi']$ norm-measurable.

We are now in a position to state the Gelfand Fundamental Theorem of
Calculus.
\begin{theorem} \label{w*FTC}
Let $\Psi:[0,T] \to X^*$ be weak* measurable and let $1\leq q \leq
\infty$ be given.  We have $[\Psi]\in W^{1,q}_{w*}(0,T;X^*)$ if and
only if there exist $[\Phi] \in L^q_{w*}(0,T;X^*)$ and $\psi_0\in X^*$
such that the mapping
\begin{equation}\label{intrepr1}
  \wt{\Psi}:[0,T] \to X^* \com{given by}
  \wt\Psi(t) := \psi_0 + \wint_0^t {\Phi}(s) \; ds
\end{equation}
satisfies $\wt{\Psi} \in [\Psi]$.  

Moreover, if $[\Psi]\in W^{1,q}_{w*}(0,T;X^*)$, then there exists an
absolutely continuous representative $\bar{\Psi} \in [\Psi]$ such that
\begin{equation}\label{intrepr2}
\bar{\Psi}(t) = \bar{\Psi}(0)+ \wint_0^t \Psi'(s) \, ds.
\end{equation}
The map $\bar{\Psi}$ belongs to $L^q(0,T; X^*)$ and has variation
$\Var_{\bar{\Psi}}$ which satisfies \eqref{w*var}.
\end{theorem}

\begin{proof}
First, suppose that $\wt{\Psi} \in [\Psi]$ satisfies \eqref{intrepr1}.
Then for every $x\in X$, $t \in [0,T]$, we have
\[
  \< \wt\Psi(t),x \> = \< \psi_0, x \> +
    \int_0^t \< \Phi(s), x \> \; ds \in \RR\,,
\]
so the function $t \mapsto\< \wt\Psi(t), x \>$ is in $W^{1,q}(0,T)$, with
weak derivative $t \mapsto\< \Phi(t), x \>$.  Thus \eqref{w*wd}
holds for any $\varphi\in C_c^\infty(0,T)$, so $[\wt\Psi]\in
W^{1,q}_{w*}(0,T;X^*)$ with $[\Phi]=[\wt\Psi']$.  

Now suppose that $[\Psi]\in W^{1,q}_{w*}(0,T;X^*)$ and let
$\hat\Psi\in[\Psi]$ and $\hat\Phi\in [\Psi']$ be norm-measurable.  We
need to find the trace $\psi_0$ of $\Psi$ at $t=0$.  For each
$x \in X$ and $\varphi\in C_c^\infty(0,T)$, we have
\[
  \int_0^T \< \hat\Psi(t),x\> \, \varphi'(t)\; dt =
   - \int_0^T \< \hat\Phi(t),x\> \, \varphi(t) \; dt\,,
\]
and \eqref{w*integrbnd} implies
$\<\hat\Psi(\cdot),x\>\in W^{1,q}(0,T)$ with weak derivative
$\<\hat\Phi(\cdot),x\>$.  This implies that for each $x\in X$, there
exists a unique absolutely continuous real-valued function
$z_x:[0,T] \to \RR$ such that
\begin{equation}
\label{aecont1}
\begin{aligned}
  z_x(t) &= \< \hat\Psi(t),x\>   \quad \text{$\lambda$-a.e. $ t \in
    [0,T]$,} \com{and} \\
  z_x(t) &= z_x(0) + \int_0^t \< \hat\Phi(s),x \> \; ds
   \com{for all} t \in [0,T]\,,
\end{aligned}
\end{equation}
where $z_x(0)$ is the trace of the map $t \to \< \hat\Psi(t), x \>$ at
$t=0$.  Note that $z_x(0)$ need not equal $\< \hat\Psi(0), x \>$, and
we must show the existence of $\psi_0 \in X^*$ such that
$z_x(0)=\< \psi_0, x\>$ for all $x \in X$.

From linearity of $\hat\Psi(\cdot)$, continuity of $z_x(\cdot)$, and
\eqref{aecont1} it follows that, for any $c_1$, $c_2 \in \RR$ and
$x_1$,~$x_2 \in X$, we have
\[
  z_{c_1x_1+c_2x_2}(t) = c_1\, z_{x_1}(t) + c_2 \,z_{x_2}(t)
  \com{for all} t\in[0,T],
\]
so the map $X \ni x \to z_x(0)$ is a linear functional.
Next, for each $x \in X$, \eqref{aecont1} implies that
\[
  \int_0^T \< \hat\Psi(t), x \> \; dt =
  \int_0^T z_x(t) \; dt = z_x(0)\,T +
     \int_0^T \Big( \int_0^t \< \hat\Phi(s), x\> \;ds \Big) \; dt\,.
\]
Since $t \to \|\hat\Psi(t)\|$ and $t \to \|\hat\Phi(t)\|$ are
measurable and integrable, we conclude that
\[
  |z_x(0)| \leq \frac{1}{T} \bigg(\int_0^T \| \hat\Psi(t)\|\; dt +
           T\, \int_0^T \|\hat\Phi(t)\| \; dt \bigg)\, \|x\|\,.
\]
Thus $X \ni x \to z_x(0)$ is also bounded, and hence there is some
$\psi_0 \in X^*$ such that $z_x(0)=\<\psi_0,x \>$ for all $x \in X$. 

It follows from \eqref{aecont1} and Theorem~\ref{w*intprop} that the
map $\bar\Psi:[0,T]\to X^*$, defined by 
\[
  \bar\Psi(t) = \psi_0 + \wint_0^t \hat\Phi(s)\;ds
\]
is absolutely continuous with G-weak derivative $\hat\Phi$, and so
$[\bar{\Psi}] \in W^{1,q}_{w*}(0,T;X^*)$.  Moreover, since
$\<\bar{\Psi}(t),x\>=z_x(t)$, we have $\bar{\Psi} \in [\Psi]$.  The
rest of the theorem follows directly from Theorem~\ref{w*intprop}.
\end{proof}

To simplify notation in the sequel, we adopt the convention that,
unless otherwise specified, we identify $[\Phi]\in L^q_{w*}(0,T;X^*)$
with a norm-measurable representative $\hat\Phi$, and if
$[\Psi]\in W^{1,q}_{w*}(0,T;X^*)$, we identify $[\Psi]$ with its
absolutely continuous representative $\bar\Psi$ and $[\Psi']$ with a
norm-measurable representative $\hat\Psi'$.

\begin{corollary}
\label{w*sobregty}
Suppose that $u \in W^{1,q}_{w*}(0,T;X^*)$ with continuous
representative $\bar u$.  If there exists a strongly measurable
$f:[0,T]\to X^*$ such that $f \in [u']$, then $\bar u$ is strongly
measurable and $\bar{u} \in W^{1,q}(0,T;X^*)$.
\end{corollary}

\begin{proof}
Since $f \in [u']$ and $\bar u$ is continuous, we have
\begin{equation}\label{ubar}
  \bar{u}(t) = \bar u(0) + \wint_0^t f(s)\; ds,\quad t\in[0,T].
\end{equation}
Moreover, since $f$ is strongly measurable, it is norm-measurable, and
\[
  \int_0^T \|f(t)\|^q \, dt  = \qn{[u']}_q^q < \infty,
\]
so that $f \in L^q(0,T;X^*)$, \eqref{ubar} is a Bochner integral, and
$\bar u\in W^{1,q}(0,T;X^*)$.
\end{proof}

Our next calculus theorem is an integration by parts formula.  Let
$1 \leq p \leq \infty$ satisfy $\frac1p+\frac1q=1$, and let $X$ be a
Banach space with dual $X^*$.

\begin{theorem}\label{w*intbypartthm}
Suppose that $f \in W^{1,p}(0,T;X)$ and $\Psi \in W^{1,q}_{w*}(0,T;X^*)$.
The function $\<\Psi(t), f(t) \>:[0,T] \to \RR$ is in $W^{1,1}(0,T)$,
with  weak derivative
\begin{equation}\label{w*derproduct}
  \mbox{\rm D} \< \Psi(t), f(t) \> =
  \< {\Psi}'(t), f(t) \> + \< \Psi(t), f'(t) \>
   \quad \mbox{$\lambda$-a.e.}\,,
\end{equation}
where $\Psi'$ and $f'$ are the G-weak and B-weak derivatives of $\Psi$
and $f$, respectively.  Moreover,
\begin{equation}\label{w*intbyparts}
  \int_0^T \< {\Psi}'(t), f(t) \> \, dt
  = \<\bar{\Psi}(s),\bar{f}(s)\> \Big|_0^T -
     \int_0^T \< \Psi(t), f'(t) \>\,dt\,,
\end{equation}
where $\bar{f}$ and $\bar{\Psi}$ are the continuous representatives of
$f$ and $\Psi$, respectively.
\end{theorem}

\begin{proof}
  Set $z(t)=\< \Psi(t), f(t) \>$, which is in $L^1(0,T)$ by Theorem
  \ref{dualLpX}, assume that $f$ and $\Psi$ are continuous
  representatives, and fix representatives $f'$ and norm-measurable
  $\Psi'$.  Let $E=\bigcup_{n=1}^N(a_n,b_n)$ be a disjoint union of
  intervals.  By Theorems \ref{BFTC} and \ref{w*FTC} we have
\[
\begin{aligned}
  f(b_n)-f(a_n) &= \int_{a_n}^{b_n} f'(t) \; dt  \com{and}\\
  \< {\Psi}(b_n)-{\Psi}(a_n), f(b_n)\>
  &= \int_{a_n}^{b_n} \< {\Psi'}(t), f(b_n)\> \; dt,
\end{aligned}
\]
and so, since $f$ and $\Psi$ are continuous, we get
\[
\begin{aligned}
  \sum_{n=1}^N & |z(b_n)-z(a_n)| \\
& \leq  \sum_{n=1}^N\Big|\int_{a_n}^{b_n}\< {\Psi'}(t), f(b_n)\> \; dt\Big|
      + \sum_{n=1}^N\Big|\int_{a_n}^{b_n}\< {\Psi}(a_n), f'(t)\>\;dt\Big|\\
& \leq \|f\|_{L^{\infty}(0,T;X)}  \int_{E} \|{{\Psi}}'(t) \|\, dt
      + \|{\Psi}\|_{L^{\infty}(0,T;X^*)} \int_{E} \|{f}'(t)\| \, dt.
\end{aligned}
\]
Absolute continuity the Lebesgue integral now implies that $z(t)$ is
absolutely continuous, so also $z\in W^{1,1}(0,T)$, with weak
derivative $\D z = \frac{dz}{dt}$ almost everywhere.

Next, Theorem \ref{w*LDT}, \eqref{BLDT} and absolute continuity of
$z(t)$ imply that the set
\[
\begin{aligned}
E = \bigg\{ t: \; &\frac{dz}{dt} \;\; \mbox{exists},
  \quad \lim_{h \to 0} \frac{1}{h} \int_0^T f'(s) \;ds = f'(t)\,, \com{and}
   \\ &  \frac{1}{h}\int_t^{t+h}  \< \Psi'(s)\, f(s) \> \; ds =
        \< \Psi'(t)\, f(t) \>  \,\,\bigg\}
\end{aligned}
\]
has full measure, $\lambda([0,T]\backslash E)=0$.  For each $t \in E$,
we have
\[
\begin{aligned}
\frac{dz}{dt}(t)
& = \lim_{h\to 0} \frac{1}{h}\Big( \< {\Psi}(t+h)-{\Psi}(t), f(t) \>  + \< {\Psi}(t), f(t+h)-f(t) \>\Big)\\[3pt]
& = \< \Psi'(t), f(t) \>  + \< {\Psi}(t), f'(t) \>\,,
\end{aligned}
\]
which is \eqref{w*derproduct}.  Next, by absolute continuity,
\[
  z(T) - z(0) = \int_{0}^T \D z(t) \, dt\,,
\]
which is \eqref{w*intbyparts}.
\end{proof}

\section{Application to Hyperbolic Conservation Laws}

We now apply the calculus we have developed to the Cauchy problem for
systems of hyperbolic conservation laws in one space dimension,
\begin{equation}\label{claw}
  \del_t u + \del_x F(u) = 0\,,\quad
  u(x,0) = u^0(x),
\end{equation}
where $u = u(x,t) : \RR\times\RR_+ \to \RR^n$, $x \in \RR$, and $F \in
C^1(\RR^n;\RR^n)$.

Recall that a \emph{distributional solution} of \eqref{claw} is a
measurable function $u(x,t):\RR\times [0,T) \to \RR^n$ which satisfies
\begin{equation}
\label{dist}
\begin{aligned}
  \int_0^T\int_{\RR}  \Big( u(x,t)^{\T} \del_t\varphi(x,t)
  &+ F(u(x,t))^{\T} \del_x\varphi(x,t) \Big)\; dx\;dt\\
   &+ \int_R  u^0(x)^{\T} \varphi(0,x) \; dx = 0,
\end{aligned}
\end{equation}
for every function $\varphi\in C^1_c(\RR^2,\RR^n)$; here it is
implicitly assumed that $F(u(x,t))$ is locally integrable.
A distributional solution which is locally bounded is a
\emph{weak solution}.

In defining distributional and weak solutions, there is no distinction
between the roles of time and space, with both the solution $u$ and
test functions $\varphi$ regarded as functions on $\RR\times [0,T)$.
Our approach is different, in that we wish to understand the PDE
\eqref{claw} as an evolution equation, so we will regard the solution
as a function of time taking values in a Banach space.  In particular,
we regard the (spatial) test functions as elements of a Banach space
$X$, which contains $C_c^{\infty}$ as a dense subspace, and, since
$u(t)$ acts linearly on these test functions, it has values in $X^*$,
and we regard $u \in L^q_{w*}(0,T;X^*)$; throughout this section, we
fix the constant $1\le q\le\infty$.  Taking this point of view, we use
the convention that $u(t) \in X^*$ stands for $u(\cdot, t)$.

The critical issue is to make sense of the nonlinear flux $F(u)$ and
its derivative in the space $X^*$.  Generally the flux $F$ is given as
a function of the conserved variables $u$, regarded as a pointwise
field.  In our case we treat $u(t)\in X^*$ as a field, and we
similarly regard $F(u(t))\in X^*$ as a field in the same way,
being defined via composition.  We refer to the corresponding map
$F:X^*\to X^*$ as a \emph{flux mapping}.  We now assume that $\Omega
\subseteq \RR$ is open and that the space $X$ of test functions
contains $C_c^\infty(\Omega)$ as a dense subset.

We say that $f\in X^*$ has an $X^*$-valued distributional derivative,
written $\DD_xf \in X^*$, if, for all $\phi\in C_c^\infty(\Omega)\subset X$,
we have
\[
  \big|\<f,\phi'\>\big| \le C\,\|\phi\|_X,
\]
and in this case we define $\DD_xf$ by
\[
  \< \DD_xf,\phi \> := - \<f,\phi'\>.
\]
Note that this is not the classical distributional derivative, because
we are requiring that $\DD_xf$ be a bounded operator on $X$.

\begin{definition}
\label{def:weakstar}
Let $F: X^* \to X^*$ be a flux mapping and $u^0 \in X^*$.  The
function $u \in W^{1,q}_{w*}(0,T; X^*)$
 is called a \emph{weak* solution} of the system \eqref{claw} if
\begin{equation}
\label{wstarode}
  u'(t) + \DD_x F(u(t)) = 0 \com{in} X^*
\end{equation}
for a.e.~$t\in[0,T]$, and such that $\bar u(0)=u^0$ in $X^*$, where
$\bar{u}$ is the time continuous representative of $u$.  Equivalently,
we require
\begin{equation}
\label{wssol}
  \bar{u}(t)-\bar{u}(s) = 
  \wint_{s}^t \DD_x F(u(\tau))\; d\tau \quad \text{in $X^*$}\,\quad
  \text{for} \quad t,\ s \in [0,T].
\end{equation}
\end{definition}

We note that this is a general definition which depends on the choice
of the space $X$ of test functions as well as the growth rate $q$.

\subsection{Classical Solutions}

As a first example, recall that if the system is symmetrizable (say if
$F$ is a gradient), then well-known energy estimates yield finite time
existence of classical solutions in the space $H^s(\RR)$ for
$s>3/2$~\cite{Majda}.  Moreover, since they are differentiable, these
satisfy the appropriate integration by parts formulae.  Thus, if we
choose
\[
  X = H^{-s}(\RR),  \com{so that}  X^* = H^s(\RR),
\]
then these classical solutions can be regarded as weak* solutions on
$[0,T)$, where $T$ is the blowup time of classical solutions.

To clarify this example, we consider the scalar Burgers' equation,
\[
  u_t + (u^2/2)_x = 0, \quad u(x,0) = u^0(x),
\]
for $u\in\RR$, which is well known.  Solving by characteristics, we
have
\begin{equation}
  u(x,t) = u^0(x^0), \com{where}
  x - x^0 = u^0(x^0)\,t,
\label{char}
\end{equation}
$x^0$ being implicitly determined as $x^0=x^0(x,t)$.  Differentiating
the equation, we see that the spatial derivative
$v(x,t)=u_x(x,t)$ satisfies
\[
  v_t + u\,v_x + v^2 = 0,
\]
and solving along characteristics, we find, after simplification,
\[
  v(x,t) = \frac{v^0(x^0)}{1+t\,v^0(x^0)},
\]
valid as long as $t < -1/v^0(x^0)$; here the blowup time is
\begin{equation}
  t_b = \inf\big\{-1/v^0(x^0)\ :\ v^0(x^0) < 0 \big\}.
\label{tb}
\end{equation}
Using the quasilinear equation $u_t + u\,u_x=0$ directly, we also see
that 
\[
  u'(t) = u_t(\cdot,t) \com{with}
  u_t(x,t) = \frac{-u^0(x^0)\,v^0(x^0)}{1+t\,v^0(x^0)}.
\]
Next, differentiating \eqref{char}, we get
\[
  \frac{\del x^0}{\del x} = \frac{1}{1+t\,v^0(x^0)} \com{and}
  \frac{\del x^0}{\del t} = \frac{-u^0(x^0)}{1+t\,v^0(x^0)}.
\]
It follows that the homogeneous $\dot H^1$ norm of the
solution is
\[
  \|u(t)\|_{\dot H^1}^2 = \int v(x,t)^2\;dx
     = \int\frac{v^0(x^0)^2}{(1+t\,v^0(x^0))^3}\;dx^0\,,
\]
where we have used $dx = (1 + t\,v^0(x^0))\,dx^0$ from \eqref{char} to
change variables.  Now using \eqref{tb}, we can write
\[
  \|u(t)\|_{\dot H^1} = O(1)\,(1-t/t_b)^{-3/2},
\]
which gives the rate of blowup of $H^1$ norm.  Thus for any $q>1$, we
have $u \in L^{q}(0,\tau;H^1)$ if and only if $\tau < t_b$, the blowup
time.  Following similar steps, we can calculate higher $H^s$ norms
and blowup rates, as desired.

This example illustrates that when the class of test functions is
large, in this case $H^{-s}$, then the solution is in a smaller space,
$H^s$, in which we have only local existence.  At the other extreme, a
very small space such as $H^s$ would yield a much larger class of
solutions.  A more realistic example is obtained by using $X = C_0$ as
the space of spatial test functions, which would lead to solutions in
the space $X^*=M$ of time varying Radon measures.  This point of view
is taken in \cite{MY2}, where the notion of weak* solution allows for
vacuums in gas dynamics and fractures in elasticity.  These solutions
are not weak solutions as vacuums and fractures are represented by
delta measures; the flux mapping $F$, as well as the entropy and
entropy flux, must be extended to measures (see \cite{MY2}).  In
general, our goal is to find the right space which would lead to
global existence and stability of solutions: in this paper, we
primarily focus on weak* solutions $u \in X^* \bigcap BV^n$, with
$X=C_0$ and $X^*=M^n$.

\subsection{BV Weak* Solutions}

It is well known that in one space dimension, the correct setting for
systems of hyperbolic conservation laws is the space $BV$ of functions
of bounded variation~\cite{G,Bressan,Dafermos}.  Recall that if
$w\in BV_{loc}$ then $\DD w\in M_{loc}$, which is the dual of $C_0$.
We will thus take
\[
  X = C_0(\RR)^n, \com{so that}  X^* = M_{loc}(\RR)^n,
\]
and use $M_{loc}$ to define weak* solutions.  However, we wish to
restrict our solutions further so that the ODE remains consistent.  We
accomplish this as follows: if
$\Phi \in W^{1,q}_{w*}\big(0,T;X^*\big)$, and in addition, $\Phi$ has
values in some $Y\subset X^*$, then we write
$\Phi \in W^{1,q}_{w*}\big(0,T;Y,X^*\big)$, that is we set
\[
  W^{1,q}_{w*}\big(0,T;Y,X^*\big)
   = \Big\{ \Phi \in W^{1,q}_{w*}(0,T;X^*)\;:\;
       \Phi(t) \in Y,\  t\in[0,T] \Big\}.
\]
Note that we do not assume that $Y$ is a subspace of $X^*$, because we
use the topology of $X^*$ throughout.  We also use the notation $T^-$
to mean up to but not including $T$, so that
\[
  L^q_{w*}(0,T^-; X^*) := \{\phi:[0,T) \to X^* \; : \;
    \phi \in L^{q}_{w*}(0,\tau;X^*)\ \forall\ 0<\tau<T\},
\]
and similarly for $W^{1,q}_{w*}(0,T^-; X^*)$.

For one-dimensional systems of conservation laws, for which $BV$
solutions are appropriate, we take $X=(C_0)^n$, $X^* = M^n_{loc}$, and
$Y=BV_{loc}^n$.

\begin{definition}\label{bvw*def}
Suppose that $u^0 \in M_{loc}(\RR; \RR^n)$.  The function
\[
  u  \in  W^{1,q}_{w*}\big( 0,T^{-};  BV_{loc}^n, M^n_{loc}\big)\,
\]
is called a \emph{BV weak* solution} to the Cauchy problem
\eqref{claw} if for a.e.~$t\in (0,T)$
\begin{equation}\label{w*claw}
  u'(t) + \DD_x F(u(t)) =0  \com{in}  M^n_{loc}
\end{equation}
and such that $\bar u(0)=u^0$, where $\bar{u}$ is the
time continuous representative of $u$ in $M^n_{loc}$.
\end{definition}

In \eqref{w*claw}, $u'$ is the G-weak derivative of $u$ in the space
$W_{w*}^{1,q}(0,T^-; M^n_{loc})$, and $\DD_x F(u(t))$ is the
distributional derivative of the function $x \to F(u(x,t))$.  Since
$u(t)\in BV_{loc}^n$ for a.e.~$t$, and $F$ is $C^1$ and so Lipschitz,
for these $t$ we also have $F(u(\cdot,t))\in BV_{loc}^n$, and
Lemma~\ref{DBV} implies that $\DD_x F(u(t)) \in M^n_{loc}$, so we can
understand the equation in $M^n_{loc}$.  Moreover, in view of
separability of $C_0$, each element of $L^q_{w*}(0,T^-;M^n_{loc})$ is
norm-measurable, and \eqref{w*claw} is equivalent to requiring that
$[u']+[\D_x F(u)]=0$ as equivalence classes of
$L^{q}_{w*}(0,T;M^n(\Omega))$.  As usual, this means
\[
  \bar{u}(t)-\bar{u}(s) = 
  \wint_{s}^t \DD_x F(u(\tau))\; d\tau \in M_{loc}
  \com{for} t,\ s \in [0,T].  
\]

We note that our definition of weak* solutions implicitly requires
some regularity, namely that $u$ is absolutely continuous as a
function $t\to u(t)\in M^n_{loc}$ with values in $BV^n_{loc}$.  Also,
recall that we require only that $u'$ (and thus $\DD_x F(u)$) be
weak*-measurable, rather than strongly measurable, and the parameter
$q$ allows a range of growth rates of weak* solutions.  In this sense,
our weak* solutions are different from weak solutions, which are
locally bounded locally integrable functions $u(x,t)$ satisfying
\eqref{dist}.  Although the two notions of solution are different, we
show they yield the same solutions in the most important case.

\begin{theorem}\label{w*wisweak}
Suppose that $u \in W^{1,q}_{w*}\big( 0,T^-; BV_{loc}^n, M^n_{loc}\big)$
is a BV weak* solution to the Cauchy problem \eqref{claw}, with
continuous representative $\bar u$.  Then $\bar u$ is H\"older
continuous as a function into $L^1_{loc}(\RR;\RR^n)$, that is,
\begin{equation}
  \label{growth}
   \bar{u} \in C^{0, 1-1/q}(0,T^-;L^1_{loc}) \com{or}
   \bar{u} \in {Lip}(0,T^-;L^1_{loc}),
\end{equation}
for $1\le q<\infty$ or $q=\infty$, respectively.  The function
$\bar{u}(x,t)$ is a distributional solution of the Cauchy problem
\eqref{claw}.   In particular, if $u$ is locally bounded,
$u \in L^{\infty}_{w*}\big( 0,T^-; L_{loc}^{\infty}(\RR;\RR^n)\big)$,
then $\bar{u}(x,t)$ is a weak solution to the Cauchy problem \eqref{claw}.
\end{theorem}

\begin{proof}
According to Theorem \ref{w11w*equiv} below,
$\bar{u}\in\wt{W}^{1,q}(0,T^{-};M^n_{loc})$, which implies that for
each open $\Omega \subset\subset \RR$, there exists $\phi\in L^q$ such
that
\[
  \|\bar{u}(t)-\bar{u}(s)\|_{M^n(\Omega)} \le
  \int_s^t\phi(\tau)\;d\tau\,,\quad t,\;s \in [0,\tau].
\]
Since $\phi\in L^q$, we apply H\"older's inequality, together with
\eqref{embed}, to get
\[
   \|\bar{u}(t)-\bar{u}(s)\|_{L^1(\Omega)} =
   \|\bar{u}(t)-\bar{u}(s)\|_{M^n(\Omega)} \le
   \|\phi\|_q\;|t- s|^{1-1/q}\,,
\]
which implies \eqref{growth}.

To show that $\bar u$ is a distributional solution, assume we are
given a test function
$\psi \in C^{1}_{c}(\RR\times (-\infty, T);\RR^n)$.  Restricting
$\psi$ to the time interval $[0,T)$, we have
$\psi\in W^{1,\infty}(0,\tau;C_0\big((a,b);\RR^n\big))$ for some
$0<\tau<T$ and finite open interval $(a,b)\in \RR$.  From
\eqref{w*claw} it follows that
$u' = -\DD _x F(u) \in L^{1}_{w*}(0,\tau;M^n(a,b))$, and by Theorem
\ref{w*intbypartthm}, the functions $t \to \< \DD_x F(u(t)),\psi(t)\>$
and $t \to \< u(t), \psi(t) \>$ are Lebesgue integrable and
\[
\begin{aligned}
  \int_0^{\tau} \< \DD_x F(u(t)),\psi(t)\> \; dt
  &= -\int_0^{\tau} \< u'(t),\psi(t) \>  \; dt \\
  &=  \<\bar{u}(0),\psi(0)\> + \int_0^{\tau} \< u(t),\psi'(t)\>\;dt\,.
\end{aligned}
\]

Next, for each $t \in [0,T)$ we have $u(t) \in BV_{loc}^n$, and since
$F$ is $C^1$, we have also $F(u(t)) \in BV_{loc}^n$, so that
$\DD_xF(u(t)) \in M^n_{loc}$ and
\[
\begin{aligned}
  \< \DD_x F(u(t)) , \psi(t) \>
  &= \int_{\RR} \psi(x, t) \, \DD_x (F(u))(dx,t) \\
  &= - \int_{\RR} F(u(x,t))^{\T} \del_x\psi(x,t)\; dx\,.
\end{aligned}
\]
Similarly, since $u(t) \in M^n_{loc} \bigcap BV_{loc}^n \subset
L_{loc}^1(\RR;\RR^n)$,
\[
\begin{aligned}
\< u(t) , \psi'(t) \>  &= \int_{\RR} \del_t \psi(x,t) \cdot u(dx,t) =
\int_{\RR}  u(x,t)^{\T} \del_t \psi(x,t) \; dx, \com{and} \\
\< \bar{u}(0) , \psi(0) \>  &= \int_{\RR} \psi(x,0) \cdot
\bar{u}(dx,0) = \int_{\RR} \bar{u}(x,0)^{\T} \psi(x,0) \; dx\,.
\end{aligned}
\]
Combining the above relations we obtain
\[
\begin{aligned}
\int_0^T \int_{\RR} \Big( u(x,t)^{\T} \, \del_t\psi(x,t) &+
F(u(x,t))^{\T} \del_x \psi(x,t) \Big) \; dx \; dt \\ &+ \int_{\RR}
u^0(x)^{\T} \psi(x,0)  \; dx = 0\,,
\end{aligned}
\]
which is \eqref{dist}.  Finally, since $\bar u(t)$, and hence also
$F(\bar u(t))$, is continuous with values in $L^1_{loc}$, both $\bar
u(x,t)$ and $F(\bar u)(x,t)$ are locally integrable, and the proof is
complete.
\end{proof}

Having shown that any weak* solution is a distributional solution, we
now show that a weak solution with sufficiently regular growth is also
a weak* solution.  In particular, the global weak solutions generated
by Glimm's method, front tracking, and vanishing viscosity, all of
which have uniformly bounded total variation, are all weak* solutions.
It follows that the uniqueness and $L^1$-stability results of Bressan
et.al.~hold unchanged in the framework of weak* solutions.

\begin{theorem}
\label{distrisw*}
Let $u(x,t)$ be a weak solution of the Cauchy problem \eqref{claw},
with $u(\cdot,t)\in BV_{loc}^n$ for each $t\in[0,T)$.  Suppose also
that for each open interval $\Omega \subset \subset \RR$, there is
some $g_{\Omega}\in L^{q}(0,T^{-})$ such that
\begin{equation}
  \label{uvarbnd}
  \Var(u(\cdot,t); \, \Omega)  \leq g_{\Omega}(t)\,,
   \quad \text{a.e.  $t \in (0,T)$.}
\end{equation}
Then $u \in W^{1,q}_{w*}(0,T^-; BV_{loc}^n, M^n_{loc})$, and $u$ is a
BV weak* solution to the Cauchy problem.
\end{theorem}

\begin{proof}
Let $ I =(a,b) \subset\subset (0,T)$ and $\Omega \subset \RR$ be
finite intervals.  Since $u(x,t)$ is locally bounded and integrable,
for each $\alpha \in C^1_0(\Omega)$, the map $t \to
\<u(t),\alpha\> := \int_{\Omega} u(x,t)\,\alpha(x)\;dx$ is measurable
and integrable over $I$.  Also,
\[
  \esssup_{t \in I} \|u(\cdot,t)\|_{M^n(\Omega)}
   = \esssup_{t \in I} \int_{\Omega} |u(x,t)| \; dx
   \leq \, \lambda(\Omega)\;\|u\|_{L^{\infty}(I \times \Omega)},
\]
so that $u \in L^{\infty}_{w*}(a,b; M(\Omega))$.

Next, for each $t$, we have $u(t):=u(\cdot,t) \in BV_{loc}^n$,
and so since $F\in C^1$, also $F(u(t)) \in BV_{loc}^n$, and in
particular, $x \to F(u(x,t))$ is integrable over $\Omega$ for each $t
\in I$.  Lemma~\ref{DBV} now implies that $\DD_x F(u (t)) \in
M^n_{loc}$ for each $t\in I$, and
\[
  t \to \<\DD_x F(u(t)),\alpha\> = -\int_{\RR} F(u(x,t))\,\alpha'(x)\;dx
\]
is measurable and integrable on $I$.  Moreover, since $F$ is locally
Lipschitz and $u$ is locally bounded, our assumption \eqref{uvarbnd}
implies that there is some $\wt g_\Omega\in L^q(0,T^-)$ such that
\[
  \|\DD_x F(u(t))\|_{M^n(\Omega)}
  = \Var\big(F(u(\cdot,t)); \Omega\big)
  \le \wt g_\Omega(t),
  \quad \text{a.e.  $t \in (0,T)$,}
\]
and so also
\[
  \big|\<\DD_x F(u (t)), \alpha \>\big|
  \leq \|\DD_xF(u(t))\|_{M^n(\Omega)}\,\|\alpha\|_{L^{\infty}(\Omega)}
  \leq {\wt{g}}_{\Omega}(t)\;\|\alpha\|_{L^{\infty}(\Omega)}\,.
\]
Since $\alpha$ is arbitrary and $g_{\Omega}\in L^q(0,T^-)$,
Theorem~\ref{w*abscont} below implies that
\[
  \DD_x F(u(t)) \in L^{q}_{w*}(0,T^-;  M^n(\Omega))\,.
\]

Now let $\psi \in C^1_c(I)$ and $\alpha \in C^1_0(\Omega)$, and set
$\varphi(x,t)=\psi(t)\alpha(x)$.  Using \eqref{dist}, we get
\[
  \int_I \Big(\int_\Omega  u(x,t)^{\T} \alpha(x)\;dx \Big)
         \, \psi'(t)\; dt
  + \int_I \Big( \int_{\Omega} F(u(x,t))^{\T} \alpha'(x)\;dx \Big)
         \, \psi(t) \; dt= 0,
\]
or equivalently
\[
  \int_I \< u(t), \alpha \>\, \psi'(t) \; dt =
  \int_I \< \DD_x F(u(t)), \alpha \>  \, \psi(t) \; dt.
\]
Since $\psi$ is arbitrary, it follows that the G-weak derivative $u'$
exists and
\[
  \<u'(t),\alpha\> = -\<\DD_x F(u(t)), \alpha\>
  \quad \text{a.e.~$t \in I$.}
\]
Since $\alpha$ and $I=(a,b) \subset\subset (0,T)$ are arbitrary we
conclude that
\[
  u'+\DD_x (F(u(t))) = 0 \quad
  \text{in $M^n(\Omega)$, a.e.~$t \in (0,T)$}
\]
and $u \in W^{1,q}_{w*}(0,T^-; M^n(\Omega))$.  Moreover, \eqref{uvarbnd}
and Theorem~\ref{w*abscont} below yield $u\in L^q_{w*}(0,T^-;BV^n(\Omega))$,
so also $u \in W^{1,q}_{w*}(0,T^-;BV^n(\Omega), M^n(\Omega))$, and
$u$ is a BV weak* solution.

It remains to show that the initial data is taken on in an appropriate
sense.  Let $\bar{u}(t)$ be the continuous representative of the weak
solution $u(\cdot,t)\in M^n(\Omega)$.  We have shown that this is a
weak* solution, with initial data $\bar u(0)\in BV^n_{loc}$.  By
Theorem~\ref{w*wisweak}, this is a continuous distributional solution.
Moreover, $u(t) = \bar u(t)$ a.e.~$t$ in $M^n(\Omega)$, and so also in
$L^1_{loc}(\Omega)$, and this in turn implies $u(x,t) = \bar
u(x,t)$ almost surely as functions of both space and time.  Since both
$u$ and $\bar u$ are distributional solutions, for
any $\varphi\in C^1_c\big(\Omega \times (-\infty,\tau)\big)$, use of
\eqref{dist} then yields
\[
  \int_{\RR} \bar{u}(x,0)^{\T}\varphi(0,x)\; dx
  = \int_{\RR} u^0(x)^{\T}\varphi(0,x)\; dx,
\]
so that $\bar{u}(0,x) = u^0(x)$ a.e.~$x$, and thus $\bar u(0) = u^0
\in M^n_{loc}$.
\end{proof}

\subsection{Shock Waves and Entropy}

It is well known that discontinuities in weak solutions satisfy the
Rankine-Hugoniot jump conditions, and that entropy conditions are
required to select the unique, physically relevant solution when
discontinuities are present.  Here we restate these conditions from
the point of view of weak* solutions.  In particular, we note the ease
with which the Rankine-Hugoniot conditions are derived.

In deriving the shock conditions, we wish to understand the local
pointwise structure of BV weak* solutions.  To do so, we make some
reasonable simplifying assumptions.  For a given $\Omega\subset\RR$
and $t\in[0,T)$, we assume $u(t)\in SBV(\Omega)$, so that $u$ has no
singular continuous part, and we can write
\begin{equation}
  u(x,t) = \sum_{j} v_j(t)\,H(x-x_j(t)) + u_c(x,t),
  \quad\text{a.e.~$x\in\Omega$,}
\label{utx}
\end{equation}
where $u_c$ is the absolutely continuous part of $u$, and there are
jumps of size $v_j(t)$ located at points $x_j(t)\in\mc
J\subset\Omega$.  The jump $v_j(t)$ is
\[
  v_j(t) = u(x_j(t)+,t) - u(x_j(t)-,t)\,,
\]
while $u_c(\cdot,t)\in W^{1,1}(\RR)$ with $\del_xu=\del_xu_c$ almost
surely, where we have used $\del_x$ to denote the pointwise (partial)
derivative.  It follows that $F(u)$ has the same form, namely
\begin{equation}
  F(u(x,t)) =  \sum_{j} g_j(t)\,H(x-x_j(t)) + f_c(x,t),
  \quad\text{a.e.~$x\in\Omega$,}
\label{Futx}
\end{equation}
where $g_j(t)$ is the jump in $F(u)$ at $x_j(t)$,
\[
  g_j(t) = F(u(x_j(t)+,t)) - F(u(x_j(t)-,t))\,,
\]
and $f_c$ is given by
\[
  f_c(x,t) =
    F\Big(u_c(x,t)+\sum_{x_j\le x}v_j(t)\Big) - \sum_{x_j\le x}g_j(t),
\]
so also $f_c(\cdot,t)\in W^{1,1}(\RR)$.  From \eqref{Futx}, we
calculate the distributional derivative
\begin{equation}
  \label{DDFut}
  \DD_x F(u) = \sum_{j} g_j(t)\,\delta_{x_j(t)}
      + \del_xf_c(\cdot,t) \in M^n(\Omega),
\end{equation}
where $\del_xf_c(\cdot,t)$ is defined $\lambda$-a.e., and in fact, for
almost all $x\ne x_j(t)$,
\[
  \del_xf_c(x,t) = DF(u)\,\del_xu_c(x,t)
              = DF(u)\,\del_xu(x,t)\in L^1(\RR),
\]
where $DF(u)$ is the derivative of $F:\RR^n\to\RR^n$.

We now make the further assumption that $u_c(x,\cdot)$, $v_j$ and
$x_j$ are absolutely continuous functions of $t$, a.e~$x\in\Omega$.
Thus for a full measure set of $t$, say $t\in E$, the weak derivative
$\DD_tu$ can be calculated as a measure, and we have
\begin{equation}
  \label{DDut}
\begin{aligned}
  \DD_t u &= \sum_{j} v_j(t)\,\frac{d(-x_j)}{dt}\,\delta_{x_j(t)}\\
     &\qquad{} + \sum_{j} \frac{dv_j}{dt}\,H(x-x_j(t)) + \del_t u_c(x,t)
  \in M^n(\Omega).
\end{aligned}
\end{equation}
Since $u$ is assumed to be a weak* solution, it follows that its
G-weak derivative $u'\in M^n(\Omega)$, and $u' = \DD_t u$ whenever
$\DD_t u$ exists as a measure, and therefore that for $t\in E$,
\[
  \DD_t u + \DD_x F(u) = 0 \in M^n(\Omega).
\]
For two measures to be equal their atomic parts must coincide, and
comparing \eqref{DDFut} and \eqref{DDut} yields
$v_j(t)\,\frac{dx_j}{dt}=g_j(t)$, which is the Rankine-Hugoniot
condition,
\begin{equation}
  \label{RH}
   F(u(x_j(t)+,t)) - F(u(x_j(t)-,t)) =
    x_j'(t)\,\big( u(x_j(t)+,t) - u(x_j(t)-,t) \big),
\end{equation}
for each $j\in\mc J$.
Moreover, away from the jump set, the absolutely continuous parts of
the measures must agree a.e., so that for a.e.~$x\ne x_j(t)$,
\[
  \sum_{x_j<x} \frac{dv_j}{dt} + \del_t u_c(x,t) + \del_xf_c(x,t) = 0,
\]
and using \eqref{utx}, this yields
\begin{equation}
  \label{QL}
  \del_t u(x,t) + DF(u)\,\del_xu(x,t) = 0,
  \quad\text{a.e.~$x\ne x_j$,}
\end{equation}
so that, as is well known, the absolutely continuous part of the
solution satisfies the quasilinear form of the equation almost
everywhere.  Note that the jump conditions \eqref{RH} are trivially
satisfied at any point of continuity of $u$.

Having checked that discontinuities in weak* solutions satisfy the
Rankine-Hugoniot condition \eqref{RH}, we must now give an entropy
selection criterion which will ensure uniqueness of solutions with
shocks.  For abstract systems that are genuinely nonlinear and/or
linearly degenerate, there are two such conditions, namely the Lax
entropy condition and that obtained from a convex entropy/flux pair.
Here we give the appropriate statement of these entropy conditions
for weak* solutions.

Recall that the Lax entropy condition holds for a single isolated jump
associated with a specific wave family.  We can rewrite \eqref{RH} for
a single isolated jump located at $x=\xi(t)$ as
\begin{equation}
  \label{RH1}
  F(u_+) - F(u_-) = \xi'(t)\;(u_+ - u_-),  \com{where}
  u_\pm = u(\xi(t)\pm,t),
\end{equation}
which is an eigenvalue problem, with $n$ different solutions, one for
each wave family.  The wave families are in turn distinguished by the
associated nonlinear wave speed $\lambda_k$, the $k$-th eigenvalue of
the matrix $DF(u)$.  Lax's entropy criterion for a $k$-shock is then
usually written as
\begin{equation}
  \label{Lax}
  \lambda_k(u_-) > \xi'(t) > \lambda_k(u_+),
\end{equation}
so the associated $k$-th characteristics impinge on both sides of the
shock.  We continue to use this condition for isolated shocks in weak*
solutions, for which the left and right limits are well-defined functions.

On the other hand, if we are given a $C^1$ entropy/flux pair
$(\eta,q)$, we require that the map
\[
  t \mapsto \eta(u(\cdot,t)) \in W^{1,1}_{w*}(0,T;M_{loc}),
\]
and the usual entropy condition holds, namely
\begin{equation}
  \label{entropy}
  \mu_\eta := \eta(u)' + \DD_x q(u) \le 0 \quad\text{in $M$.}
\end{equation}
This can be interpreted as a measure, since every signed
distribution is representable as a signed measure.  Note that since
$\eta$ and $q$ are $C^1$ functions, use of \eqref{QL} at points of
absolute continuity of $u(x,t)$ means that the measure $\mu_\eta$
vanishes at those points, so that $\mu_\eta$ is supported on the
jump set of the solution.  We note that the requirement that $\eta\in
W^{1,1}_{w*}(0,T;M_{loc})$ is stronger than the usual assumption that
the sum $\eta_t+q_x = \text{div}(\eta,q) \le 0$ in $M$, because we require
each term to be a measure separately. 

\subsection{The Riemann Problem}

As an illustrative example, we describe the well-known solution of the
Riemann problem as a $BV$ weak* solution.  Recall that the solution of the
Riemann problem consists of constant states separated by $n$ centered
elementary waves, these being (centered) shocks or rarefactions for
genuinely nonlinear families, and (centered) contact discontinuities
for linearly degenerate families.

Recall that a centered $k$-shock consists of a discontinuity between
two constant states $u_- = u_{k-1}$ and $u_+ = u_k$ satisfying
\eqref{RH1} and \eqref{Lax}, which has constant shock speed
$\xi'(t) = \sigma_k(u_{k-1},u_k)$, so that the position of the shock
at time $t$ is given by $\xi(t) = \sigma_k(u_{k-1},u_k)\,t$.  The
centered $k$-shock can thus be written as
\begin{equation}
  \label{kshock}
  u(x,t) = u_{k-1} + (u_k-u_{k-1})\,H(x-\sigma_k(u_{k-1},u_k)\,t).
\end{equation}

A centered rarefaction wave is a solution of the form
$u(x,t)=w(x/t)$, which satisfies \eqref{QL}, which reduces to
\[
  -\eps\,\frac{dw}{d\eps} + DF(w(\eps))\,\frac{dw}{d\eps} = 0,
  \quad \eps = x/t.
\]
This is an eigenvalue problem, and the centered $k$-rarefaction wave
is $u(x,t)=w_k(\eps)$, corresponding to the $k$-th eigenpair,
\begin{equation}
  \frac{dw_k}{d\eps}=r_k(\eps) \com{and}
  \eps = \frac xt = \lambda_k(w_k(\eps)),
  \label{kwave}
\end{equation}
and such that $\lambda_k(w_k(\eps))$ is monotone increasing.  The
$k$-rarefaction between states $u_{k-1}$ and $u_k$ is then described by
\begin{equation}
  u(x,t) =
  \begin{cases}
    u_{k-1}, & x/t \le \lambda_k(u_{k-1}), \\
    w_k(\eps), &
       \lambda_k(u_{k-1}) \le x/t = \eps \le \lambda_k(u_k), \\
    u_k, & \lambda_k(u_k) \le x/t.
  \end{cases}
  \label{kRf}
\end{equation}

If the $k$-th family is linearly degenerate, that is
$r_k\cdot\nabla\lambda_k\equiv0$, then the wave speed does not change
across a $k$-wave, and both \eqref{kshock} and \eqref{kRf} degenerate
and coincide, and the $k$-contact discontinuity is given by
\begin{equation}
\begin{aligned}
  u(x,t) &= u_{k-1} + (u_k-u_{k-1})\,H(x-\lambda_k(u_k)\,t),
  \\ \lambda_k(u_k) &= \lambda_k(u_{k-1}).
\end{aligned}
  \label{kcontact}
\end{equation}

The general Riemann problem with data $u^0(x) = u_L + (u_R-u_L)\;H(x)$
is solved by finding (unique) states $\{u_\ell\}_{\ell=0}^n$, with
$u^0=u_L$ and $u_n=u_R$, such that each $u_{k-1}$ is connected to
$u_k$ by a centered $k$-wave.  We describe this as a weak* solution by
explicitly calculating the appropriate derivatives.  It suffices to
consider the waves separately.

We first consider a $k$-shock satisfying \eqref{kshock}.  Clearly
\[
  F(u(x,t)) = F(u_{k-1})
     + (F(u_k)-F(u_{k-1}))\,H(x-\sigma_k(u_{k-1},u_k)\,t),
\]
and differentiating yields
\[
\begin{aligned}
  u' &= \DD_t u = - (u_k-u_{k-1})\,\sigma_k\,\delta_{\sigma_k t}
\com{and}\\
  \DD_x F(u) &= (F(u_k)-F(u_{k-1}))\,\delta_{\sigma_k t},
\end{aligned}
\]
where $u'$, being strongly measurable in $M_{loc}^n$, is the B-weak
derivative of $u$.  The jump condition \eqref{RH1} implies
\eqref{w*claw}, and $u'$ is clearly bounded,
\[
\begin{aligned}
  \big|\< u', \alpha \>\big| &=
  |u_k-u_{k-1}|\,|\sigma_k|\,\|\alpha\|, \com{and}\\
  \|u'\|_{M^n} &= |u_k-u_{k-1}|\,|\sigma_k|,
\end{aligned}
\]
so $u \in W^{1,\infty}\big( 0,\infty^{-};  BV_{loc}^n, M^n_{loc}\big)$.
The same estimates hold for $k$-contacts.

We now describe the $k$-rarefaction similarly.  Across the wave,
namely for $\lambda_k(u_{k-1}) \le x/t \le \lambda_k(u_k)$, we have
\[
  u' = \frac{dw_k}{d\eps}\,\del_t\eps  \com{and}
  \DD_x F(u) = DF\,\frac{dw_k}{d\eps}\,\del_x\eps,
\]
where $x = \lambda_k(w_k(\eps))\,t$, and $u' = \DD_x F(u) = 0$
otherwise.  It follows that
\[
  1 = \dot\lambda_k\,t\,\del_x\eps \com{and}
  0 = \lambda_k + \dot\lambda_k\,t\,\del_t\eps,
\]
where $\dot\lambda_k:=\frac{d\lambda_k}{d\eps}
=\frac{dw_k}{d\eps}\cdot\nabla\lambda_k$, and we can calculate the
action of $u'$ on a test function $\alpha$, by
\[
  \<u',\alpha\> =
  \int_{t\lambda_k(u_{k-1})}^{t\lambda_k(u_k)} \frac{dw_k}{d\eps}\,
  \frac{-\lambda_k}{t\dot\lambda_k}\,\alpha(x)\;dx
  = -\int_{\lambda_k(u_{k-1})}^{\lambda_k(u_k)}
  \frac{dF}{d\lambda}\,\alpha(\lambda t)\;d\lambda,
\]
by change of variables, where now $\lambda$ parameterizes the integral
curve.  It follows that $u'$ is bounded,
\[
  \big|\<u',\alpha\>\big| \le
  \int_{\lambda_k(u_{k-1})}^{\lambda_k(u_k)}
  \Big|\frac{dF}{d\lambda}\Big|\;d\lambda\;\|\alpha\|, \com{so that}
  \|u'\|_{M^n} \le \Var(F),
\]
the variation being taken along the integral curve, and again the
solution satisfies $u\in
W^{1,\infty}\big( 0,\infty^{-};  BV_{loc}^n, M^n_{loc}\big)$.

Having described the simple waves separately, we now combine them into
Lax's well-known solution of the general Riemann problem.  That is, we
identify states $u_0=u_L,\ u_1,\dots,\ u_n=u_R$, such that each pair
$(u_{k-1},u_k)$ is joined by a $k$-wave, satisfying \eqref{kshock} or
\eqref{kRf} for genuinely nonlinear fields, and \eqref{kcontact} for
linearly degenerate fields.  Our assumption of strict hyperbolicity
means that for each $k=1,\dots,n$, we have
$\lambda_k(u_k)<\lambda_{k+1}(u_k)$, so the waves can be consistently
pieced together for positive times $t>0$.  This means that the
calculations above all hold locally, and for $t>0$, the (Bochner)
derivative $u'$ is simply the sum of each of the individual terms.

It remains to state the sense on which the initial data is taken on.
Recall that the solution is absolutely continuous in $M^n$, and as
$t\to0$, each individual wave converges in $M^n$ to the limit
$u_{k-1}+(u_k-u_{k-1})\,H(x)$, so the full solution converges to
$u_L+(u_R-u_L)\,H(x)$, as required.  Here, since for each fixed $t$,
$u(\cdot,t)\in BV\subset L^1_{loc}\subset M^n_{loc}$, we note that
$\|u(\cdot,t)\|_{M^n} = \|u(\cdot,t)\|_{L^1}$, so we have continuity
in $L^1$ as a function of $t$, namely $u(\cdot,t)\to u^0\in L^1$.  We
note that generally, the solution is not continuous at $t=0$ in $BV$,
but only in the larger space $L^1$.  However, the Riemann solution is
continuous as a function of $t$ in $BV$ (with constant $BV$ norm) for
all positive times $t>0$.

\appendix
\section{G-weak Differentiability}

We now show that if a function $\Psi \in L^{q}_{w*}(0,T;X^*)$ is
bounded by a sufficiently regular integral, then it is
G-weakly differentiable.  Throughout this appendix, we assume $p$ and
$q$ satisfy $\frac1p+\frac1q=1$.

\begin{theorem}\label{w*abscont}
Fix $\psi_0 \in X^*$ and suppose that $\Psi \in L^{q}_{w*}(0,T;X^*)$
is such that for each $x\in X$, there are functions $v_x \in L^1(0,T)$
such that
\begin{equation}\label{psiabscont}
  \< \Psi(t),x \>  = \< \psi_0, x \> + \int_0^t \, v_x(s) \; ds\,,
  \quad \text{a.e. $t \in [0,T]$},
\end{equation}
and suppose there is a non-negative $v\in L^q(0,T)$  such that
\begin{equation}\label{varintbnd}
|v_x(t)| \, \leq v(t) \,\|x\|\,,  \quad \text{a.e. $t \in [0,T]$}\,.
\end{equation}
Then we have
$[\Psi] \in W^{1,q}_{w*}(0,T;X^*)$  if either $1<q\le\infty$ or $q=1$
and $X$ is separable.
\end{theorem}

\begin{proof}
  Recalling that $L^{p}(0,T;X)^*\simeq L^{q}_{w*}(0,T;X^*)$, we
  construct a bounded linear functional on the space $L^{p}(0,T;X)$
  which is a G-weak derivative of $[\Psi]$.

  We begin by observing that $v_x$ is almost linear in $x$; that is,
  for fixed $x_1$, $x_2\in X$ and scalars $\alpha_1$, $\alpha_2$, we
  have
\[
  \int_0^t v_{\alpha_1 x_1+\alpha_2 x_2}(s)\;ds
  = \int_0^t \alpha_1 v_{x_1}(s)+\alpha_2v_{x_2}(s)\;ds,
    \quad \text{a.e. $t$,}
\]
and so absolute continuity of the integral implies that, for a.e.~$t$,
we have 
\[
  v_{\alpha_1\, x_1+\alpha_2\, x_2}(t) = 
  \alpha_1\, v_{x_1}(t) + \alpha_2\, v_{x_2}(t).
\]

We first assume $q>1$ and let $h(t) = \sum_{n=1}^N x_n \,\mc
X_{E_n}(t)\in X$ be a simple function.  
Define the functional $\Gamma$ on simple
functions by
\begin{equation}\label{funcdefsimple}
  \< \Gamma, h \> := \sum_{n=1}^N \int_{E_n} v_{x_n}(t) \; dt\,.
\end{equation}
It is easy to verify that $\<\Gamma,h\>$ is independent of the
representation of $h$, and
\begin{equation}
  \< \Gamma, \alpha_1 h_1 + \alpha_2 h_2 \>
   = \alpha_1 \< \Gamma, h_1 \> + \alpha_2 \< \Gamma, h_2 \>\,,
\label{gammalins}
\end{equation}
so that $\Gamma$ is linear on the subspace of simple
functions.

Recall that for simple $h(t)$ we can take the $E_n$ disjoint, in which
case
\[
  \|h\|_{L^p(0,T;X)} =
  \Big(\int_0^T \Big\|\sum_n x_n\,\mc X_{E_n}(t)\Big\|^p \; dt
  \Big)^{1/p}
  = \Big(\sum_n \|x_n\|^p \lambda(E_n)  \Big)^{1/p}\,,
\]
and for such $h$, we have
\begin{equation}\label{gammabnds}
\begin{aligned}
  \big|  \< \Gamma, h \> \big| \, & \leq
   \sum_{n=1}^{N}\|x_n\|\int_{E_n} v(t) \; dt
   = \,\int_0^T v(t) \Big(\sum_{n=1}^{N}\|x_n\|\mc X_{E_n}(t) \Big) \; dt \\
& \leq \|v\|_{L^q(0,T)}
     \Big(\sum_{i=1}^N \|x_n\|^p \lambda(E_n) \; dt \Big)^{1/p}
    = \| v \|_{L^q(0,T)} \|h\|_{L^p(0,T;X)}\,,
\end{aligned}
\end{equation}
where we have used \eqref{varintbnd} and Young's inequality.

If $f\in L^p(0,T;X)$, take a sequence $\{h_n\}_{n \geq1}$ such that
$\|f-h_n\|_{L^p(0,T;X)} \to 0$.  By \eqref{gammabnds}, the sequence
$\big\{\< \Gamma, h_n \> \big\}_{n \geq 1}$ is Cauchy and hence
converges, so we set $\< \Gamma, f \> := \lim_{n \to \infty} \<
\Gamma, h_n \>$, this limit being independent of the sequence $h_n$.
Using \eqref{gammalins}, it follows that $\Gamma$ is linear, and from
\eqref{gammabnds} it is bounded,
\[
  \big|\< \Gamma, f \>\big| \leq
  \big\|v\|_{L^q(0,T)} \lim_{n \to \infty} \|h_n\|_{L^p(0,T;X)}
  = \|v\|_{L^q(0,T)} \|f\|_{L^p(0,T;X)}\,,
\]
so that $\Gamma \in L^p(0,T;X)^*$.  By Theorem~\ref{dualLpX}, there
exists $[\Phi] \in L^q_{w*}(0,T;X^*)$ such that
\[
  \< \Gamma, f \> = \int_0^T \< {\Phi}(s), f(s) \> \; ds
   \com{for all} f \in L^p(0,T;X)\,.
\]
Define $\bar{\Psi}:[0,T] \to X^*$ by 
\[
  \bar{\Psi}(t) = \psi_0 + \wint_0^t {\Phi}(s) \; ds, \com{so that} 
  [\bar{\Psi}] \in W^{1,q}_{w*}(0,T;X^*).
\]
Then, recalling \eqref{psiabscont} and
\eqref{funcdefsimple}, we have for each $x \in X$
\[
\begin{aligned}
  \< \bar\Psi(t), x \> &
  = \<\psi_0, x\> + \int_0^T \<{\Phi}(s),x\,\mc X_{[0,t]}(s) \> \; ds
  = \< \psi_0, x \> + \< \Gamma, x\,\mc X_{[0,t]} \> \\
  &=  \< \psi_0, x \> + \int_{0}^t v_{x}(s)\; ds = \< {\Psi}(t),x
  \>\,,
 \quad \text{a.e. $t \in [0,T]$},
\end{aligned}
\]
so that $\Psi \in [\bar{\Psi}]$ and thus $\Psi \in
W^{1,q}_{w*}(0,T;X^*)$, completing the proof for $q>1$.

Now set $q=1$ and suppose that $\{x_n\}_{n \geq 1}$ is dense in $X$.
Then the set
\[
  E = \bigcup_{n \geq 1} \Big\{t \in (0,T): \,
   \lim_{h \to 0}\frac{1}{h} \int_t^{t+h} v_{x_n}(s)\;ds = v_{x_n}(t)\,,
\quad |v_{x_n}(t)| \leq v(t)\,\|x_n\| \Big\}
\]
has full measure.  For $t_0 \in E$, define $\Phi(t_0)$ by
$\<\Phi(t_0),x_n\> = v_{x_n}(t_0)$ for each $n$.  Then, as above, we
check that $\Phi(t_0)$ is a bounded linear functional on the linear
span of the set $\{x_n\}_{n \geq 1}$, with norm bounded by $v(t_0)$.
This can be extended by continuity to a linear functional on all of
$X$, with norm $\|\Phi(t_0)\|_{X^*}\le v(t_0)$.  Setting $\Phi(t_0)=0$
for $t_0\notin E$, we obtain a mapping $\Phi:[0,T] \to X^*$ such that
$\< \Phi(t), x_n\> = v_{x_n}(t)$ a.e.~$t \in [0,T]$ for each $n$, and
$\lim_{n\to\infty}\< \Phi(t), x_n\>=\<\Phi(t),x\>$ for each $t$.  This
implies that $\Phi$ is weak* measurable, and $\|\Phi(t)\| \leq v(t)$
implies that $[\Phi] \in L^1_{w*}(0,T;X^*)$.  Again define
$\bar{\Psi}:[0,T] \to X^*$ by $\bar{\Psi}=\psi_0+\wint_0^t \, \Phi(s)
\; ds$, so that $[\bar{\Psi}] \in W^{1,1}(0,T;X^*)$.

For fixed $x\in X$, the set
\[
  A_x = \Big\{ t : \<\Psi(t),z\>=\<\psi_0,z\>+\int_0^t v_z(s)\;ds\,,
  \com{for}  z=x\ \text{and each}\ x_n\Big\}
\]
has full measure.  For $t \in A_x$ and any $n$, we have
$\<\bar{\Psi}(t),x_n\>=\<\Psi(t),x_n\>$, so that
\[
\begin{aligned}
  \big|\< \bar{\Psi}(t) - \Psi(t),x\>\big| &
  = \big|\<\bar{\Psi}(t),x-x_n\>-\<\Psi(t),x-x_n\>\big| \\
  &  \leq \|\bar{\Psi}(t)\|\|x-x_n\| +\|\Psi(t)\|\|x-x_n\|,
\end{aligned}
\]
so that $\<\bar{\Psi}(t),x\> = \<\Psi(t),x\>$ for all $t \in A_x$.
Since $x\in X$ is arbitrary, $\Psi$ and
$\bar{\Psi}$ are weak*-equivalent and so
$\Psi\in[\bar\Psi] \in W^{1,1}_{w*}(0,T;X^*)$.
\end{proof}

\subsection{Relation to Brezis' Space}

Given a Banach space $X$ and $1\le p\le \infty$, Brezis defined the
space $\wt{W}^{1,p}(0,T;X)$ to be the set of all absolutely continuous
functions $\Psi:[0,T] \to X$ for which the total variation function
$t \to \Var_{\Psi}(t)$ is absolutely continuous on $[0,T]$, and such
that the scalar pointwise derivative
$\frac{d \Var_{\Psi}}{dt}\in L^p(0,T)$.  In particular, for
$\Psi\in\wt{W}^{1,p}(0,T;X)$, for all $0\le s<t\le T$, we have
\begin{equation}
  \label{varbnd}
  \|\Psi(t)-\Psi(s)\| \le
  \Var_{\Psi}(t) - \Var_{\Psi}(s) =
  \int_s^t \frac{d}{dt}\Var_{\Psi}(\tau) \; d\tau.
\end{equation}

Brezis proved that $W^{1,p}(0,T;X) \subset \wt{W}^{1,p}(0,T;X)$, and
if $X$ is reflexive, the two spaces coincide, $W^{1,p}(0,T;X) =
\wt{W}^{1,p}(0,T;X)$~\cite{Brezis73}.

We have shown in Theorems \ref{w*FTC} and \ref{w*intprop} that each
$\Psi\in W^{1,q}_{w*}(0,T;X^*)$ has an absolutely continuous
representative $\bar{\Psi}$, and that this in turn satisfies
$\bar\Psi\in\wt{W}^{1,q}(0,T;X^*)$.  We now show that the converse is
also true, that is, if $\Psi \in \wt W^{1,q}(0,T;X^*)$, then
$[\Psi]\in W^{1,q}_{w*}(0,T;X^*)$.

\begin{theorem}
\label{w11w*equiv}
For $1 < q \leq \infty$, the canonical mapping $\wt{W}^{1,q}(0,T;X^*)
\to W^{1,q}_{w*}(0,T;X)$, given by $\Psi\mapsto[\Psi]$, is injective,
onto and norm-preserving.  The same conclusion holds for $q=1$
provided $X$ is separable.
\end{theorem}

\begin{proof}
Suppose that $\Psi\in\wt{W}^{1,q}(0,T;X^*)$ and fix $x \in X$.
Then the numerical function
\[
  z_x(t) := \< \Psi(t) , x \>:[0,T] \to \RR
\]
is absolutely continuous, so its pointwise derivative
$\frac{dz_x}{dx}$ is defined $\lambda$-a.e., belongs to $L^1(0,T)$,
and
\[
  \< \Psi(t), x\> = \< \Psi(s), x\> +
   \int_s^t \frac{d z_x}{dt}(\tau) \; d\tau \com{for all}
   0\le s < t \le T\,.
\]
Now, using \eqref{varbnd},  we have for all $0 \le s < t \le T$,
\[
  \Big| \int_s^t \frac{d z_x}{dt}(\tau) \; d\tau  \Big| =
   \big| \< \Psi(t) - \Psi(s), x \> \big| \le
   \|x\| \int_s^t \frac{d}{dt} \Var_{\Psi}(\tau) \; d\tau,
\]
which in turn implies
\[
  \Big|\frac{dz_x}{dt}(t)\Big| \le
  \|x\|\,\frac{d}{dt}\Var_{\Psi}(t) \quad \text{a.e.  $t\in[0,T]$.}
\]
Recalling that the variation satisfies
$\frac{d}{dt}\Var_{\Psi}(\cdot) \in L^q(0,T)$, the result now follows
from Theorem \ref{w*abscont}.
\end{proof}

\bibliographystyle{amsplain}

\end{document}